\documentclass[11pt]{article}

\usepackage[utf8]{inputenc}
\usepackage[T2A]{fontenc}
\usepackage[russian]{babel}
\usepackage{graphicx}
\usepackage{epstopdf}
\usepackage{amsmath, amsthm, amssymb}

\graphicspath{{ITypes/}{TTypes/}{Theorem1/}{Theorem2/}{Theorem3/}}
\makeatletter
\def\th@plain{
  \thm@notefont{}
  \itshape
}
\def\th@remark{
  \thm@headfont{\itshape}
  \thm@notefont{\itshape}
  \thm@bodyfont{\normalfont}
}
\makeatother

\theoremstyle{plain}
\newtheorem{lemma}{Лемма}[section]
\newtheorem{theorem}{Теорема}
\newtheorem*{problem}{Задача}
\newtheorem*{corollary}{Следствие}
\newtheorem*{ntheorem}{Теорема}

\theoremstyle{remark}
\newtheorem*{remark}{Примечание}

\begin{document}
\centerline{\textbf{Обобщенная задача Аполлония}}
\smallskip
\centerline{Е.\,А.\,Морозов\footnote{gorg.morozov@gmail.com}}
\smallskip
\centerline{ГБОУ Лицей <<Вторая школа>>, г. Москва, 2016-2017 гг.}
\bigskip

\textbf{Abstract.} The aim of this paper is to generalize Apollonius' problem. The problem is to construct a circle that is tangent to three given circles in a plane. We find the maximum possible number of solution circles in the case of more than the three given circles. We show that if all the given circles are not tangent at the same point, then there exist at most six solutions in the case of the four given generalized circles and there exist at most four solutions in the case of the five given generalized circles. We also describe all quadruples of generalized circles with exactly six solutions.
%В работе делается попытка обобщения исследования известной задачи Аполлония о построении окружности, касающейся 3-х данных. Рассматривается вопрос о максимальном возможном числе таких окружностей, в случае если исходных окружностей больше 3-х. Доказано, что если не все исходные окружности касаются в одной точке, то в случае 4-х исходных окружностей имеется не более 6-и решений задачи Аполлония, а в случае 5-и исходных окружностей --- не более 4-х. Также дано описание всех четверок окружностей, для которых количество решений максимально.

\section{Введение}
\begin{problem}[классическая задача Аполлония]
Построить циркулем и линейкой окружность, касающуюся каждой из 3-х данных.
\end{problem}
Так сформулировал задачу известный древнегреческий геометр Аполлоний Пергский в III в. до н. э. С~тех пор было придумано много остроумных решений первоначальной задачи, однако нас будет интересовать лишь вопрос о максимально возможном числе искомых окружностей. Известно, что классическая задача имеет не более 8-и решений (см., например, \cite{bib:general}). Мы рассмотрим обобщение этого вопроса на случаи 4-х и 5-и окружностей.

Ключевую роль в наших рассуждениях играет инверсия, упрощающая конструкцию до обозримого перебора. Поэтому перед доказательством основных теорем мы подробнее рассмотрим конструкции, к которым нас приведет это преобразование (см. п.~3). Также мы используем один из результатов полного перебора всех случаев классической задачи Аполлония, описанного в \cite{bib:general}.

Среди приемов, которые мы используем для доказательства различных случаев основной теоремы будут: подсчет количества точек пересечения окружностей, их взаимное расположение на плоскости, а также факт об объединении решений классической задачи в пары, использующий алгебраическую интерпретацию окружности --- об этом подробнее в п.~4.

\section{Формулировки основных результатов}
\emph{Обобщенной окружностью (объектом)} назовем окружность, прямую или точку.

Две обобщенные окружности \emph{обобщенно касаются}, если выполнено одно из следующих условий:
\begin{enumerate}
  \item{Это две касающиеся окружности или касающиеся окружность и прямая.}
  \item{Это две параллельные прямые.}
  \item{Это точка, лежащая на окружности или прямой.}
\end{enumerate}

Так как все касания у нас будут обобщенными, то для краткости вместо «обобщенно касаются» мы будем говорить просто «касаются».

Заметим, что посредством стереографической проекции плоскость, с добавленной к ней бесконечно удаленной точкой $\infty$, можно перевести в сферу. При этом преобразовании обобщенная окружность на плоскости перейдет в окружность или точку на сфере, а касающиеся обобщенные окружности на плоскости --- в касающиеся окружности или инцидентные окружность и точку на сфере. Поэтому фактически нами рассматривается сфера и ее непустые пересечения с плоскостями. В~дальнейшем мы будем считать, что бесконечно удаленная точка на плоскости присутствует, в частности, через эту точку проходят все прямые и сама она является обобщенной окружностью. Эту точку можно вернуть на видимую плоскость с помощью инверсии, которая меняет ее местами со своим центром. Строгие определения и некоторые теоремы см. в~\cite{bib:inversion}.

Теперь можно сформулировать основные результаты.

\begin{theorem}
На плоскости даны 4 различные обобщенные окружности, которые не все касаются в одной точке. Тогда существует не более 6-и обобщенных окружностей, касающихся каждой из данных.
\end{theorem}

\begin{theorem}
На плоскости дано 5 различных обобщенных окружностей, которые не все касаются в одной точке. Тогда существует не более 4-х обобщенных окружностей, касающихся каждой из данных.
\end{theorem}

\begin{figure}[h]
\begin{minipage}[h]{0.36\linewidth}
\center{\includegraphics[width=1\linewidth]{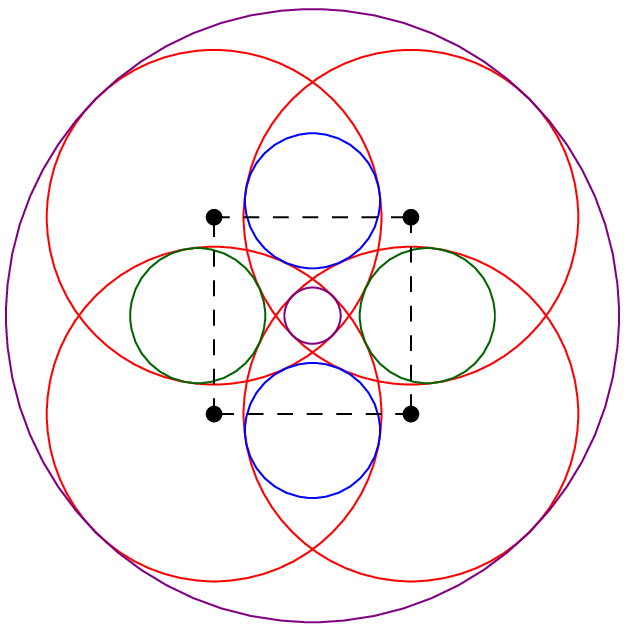}}
\end{minipage}
\hfill
\begin{minipage}[h]{0.36\linewidth}
\center{\includegraphics[width=1\linewidth]{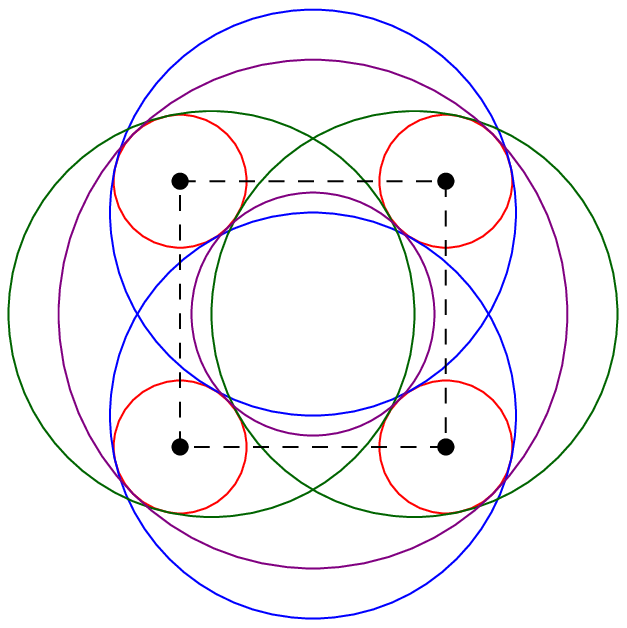}}
\end{minipage}
\caption{Примеры к теореме~1}
\label{ris:examples}
\end{figure}

Заметим, что количество окружностей, указанное в теоремах~1 и 2, действительно достигается. Соответствующие примеры для 4-х окружностей приведены на Рис.~\ref{ris:examples}. В~этих примерах исходные окружности имеют равный радиус, а их центры лежат в вершинах квадрата. Обоснование примеров практически очевидно и не нуждается в комментариях. Пример для 5-и исходных окружностей получается из данных, если рассмотреть 5 решений в качестве исходных окружностей.

В п.~7 дается полное описание четверок окружностей с 6-ю решениями.

Данные результаты, по-видимому, являются новыми.

\section{Обозначения, договоренности и используемые конструкции}
Сначала договоримся о некоторых общих обозначениях. Данные в (обобщенной) задаче Аполлония объекты назовем \emph{исходными}, а объект, касающийся всех данных, назовем \emph{решением}. Назовем конфигурацию обобщенных окружностей \emph{вырожденной}, если какие-то три из них касаются в одной точке. Вырожденная тройка, очевидно, имеет бесконечно много решений.

Мы пользуемся следующими известными свойствами инверсии (о них см., например, в \cite{bib:inversion}):
\begin{enumerate}
  \item{Пару непересекающихся окружностей или непересекающиеся окружность и прямую можно перевести инверсией в концентрические окружности.}
  \item{Пару пересекающихся окружностей или пересекающиеся окружность и прямую можно перевести инверсией в пересекающиеся прямые.}
  \item{Пару касающихся окружностей или касающиеся окружность и прямую можно перевести инверсией в параллельные прямые.}
  \item{Если точка лежит вне окружности, то существует инверсия с центром в этой точке, оставляющая окружность на месте.}
  \item{Если точка лежит внутри окружности, то существует композиция инверсии и центральной симметрии с центром в этой точке, оставляющая окружность на месте (т. н. \emph{инверсия с мнимым радиусом}).}
\end{enumerate}

Если два объекта из трех фиксированы, а третий $\alpha$ меняется, то соответствующие решения для всех трех объектов мы называем \emph{решениями для $\alpha$} (или \emph{решениями, порожденными $\alpha$}), а объект, касающийся только двух фиксированных, мы называем \emph{допустимым}.

Далее определим использующиеся в \cite{bib:general} \emph{метки Фитцджеральда} троек объектов, каждый из которых является окружностью или прямой. Такой тройке сопоставляется набор букв, характеризующий их отношения: буква $I$ используется для обозначения каждой пары пересекающихся объектов, буква $T$ --- для обозначения каждой пары касающихся объектов, а буква $S$ --- для обозначения того, что некоторая окружность или прямая разделяет два других объекта. Если все три объекта имеют общую точку, то запись заключается в квадратные скобки. Таким образом, метке $[TTT]$ соответствует вырожденная тройка объектов.

Также нам потребуются некоторые определения, связанные с тремя простейшими конструкциями из обобщенных окружностей. В~каждой из них мы фиксируем два объекта.

\textbf{1. Зафиксированы две пересекающиеся исходные прямые.}

Здесь будем называть 4 части, на которые две пересекающиеся прямые разбивают плоскость, \emph{секторами}. Сектора пронумеруем по традиции обозначения квадрантов --- римскими числами против часовой стрелки. Точку пересечения прямых обозначим через $O$.

Каждая окружность, касающаяся двух прямых, лежит в одном из секторов. Эти окружности здесь --- \emph{допустимые}. Назовем \emph{распределением} некоторого набора допустимых окружностей ненулевого радиуса по секторам запись вида $x$-$y$-$z$-$t$, где число на $i$-ом месте означает, сколько окружностей вписано в $i$-ый сектор. Разбиения \emph{изоморфны}, если они совпадают с точностью до переворотов и циклических сдвигов строки $x$-$y$-$z$-$t$, т.~е., например, разбиения 0-1-2-3 и 2-1-0-3 изоморфны.

Теперь добавим еще одну окружность или прямую $\alpha$ и рассмотрим распределение решений для $\alpha$.

\begin{figure}
\begin{minipage}[h]{0.3\linewidth}
\center{\includegraphics[width=1\linewidth]{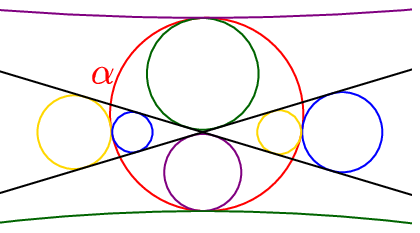}} \\ I тип
\end{minipage}
\hfill
\begin{minipage}[h]{0.3\linewidth}
\center{\includegraphics[width=1\linewidth]{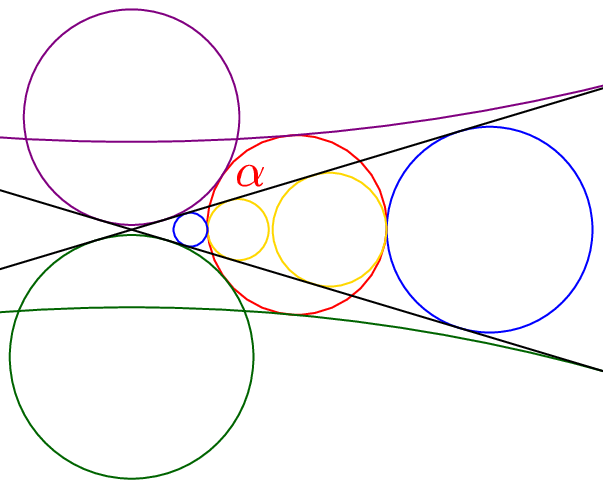}} \\ II тип
\end{minipage}
\hfill
\begin{minipage}[h]{0.3\linewidth}
\center{\includegraphics[width=1\linewidth]{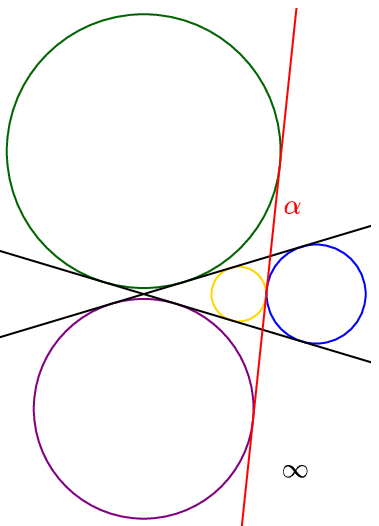}} \\ V тип, прямая
\end{minipage}
\vfill
\begin{minipage}[h]{0.3\linewidth}
\center{\includegraphics[width=1\linewidth]{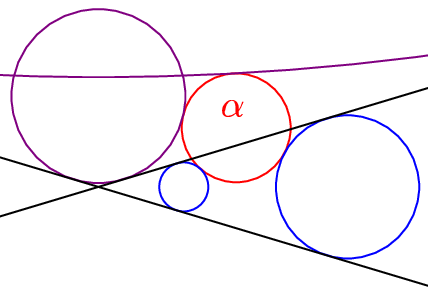}} \\ III тип
\end{minipage}
\hfill
\begin{minipage}[h]{0.3\linewidth}
\center{\includegraphics[width=1\linewidth]{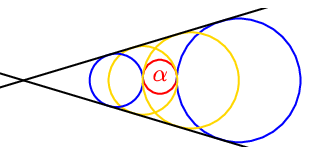}} \\ IV тип
\end{minipage}
\hfill
\begin{minipage}[h]{0.3\linewidth}
\center{\includegraphics[width=1\linewidth]{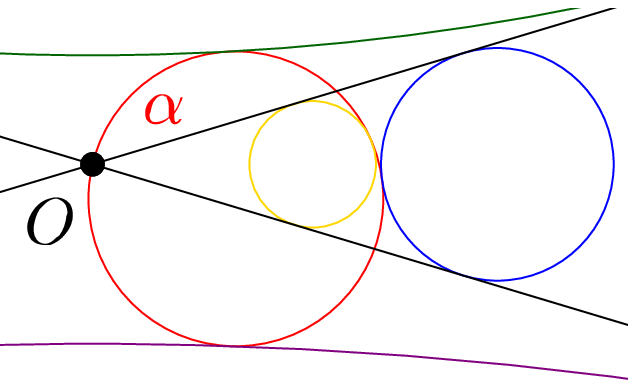}} \\ V тип, окружность
\end{minipage}
\vfill
\begin{minipage}[h]{0.3\linewidth}
\center{\includegraphics[width=1\linewidth]{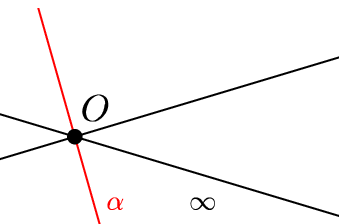}} \\ VI тип
\end{minipage}
\hfill
\begin{minipage}[h]{0.3\linewidth}
\center{\includegraphics[width=1\linewidth]{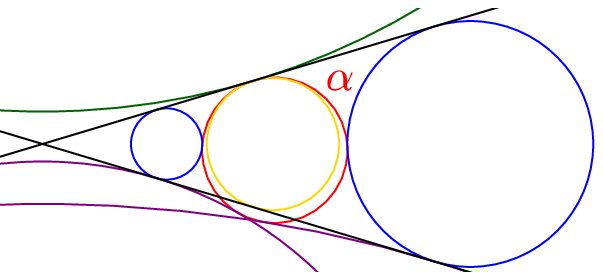}} \\ VII тип
\end{minipage}
\hfill
\begin{minipage}[h]{0.3\linewidth}
\center{\includegraphics[width=1\linewidth]{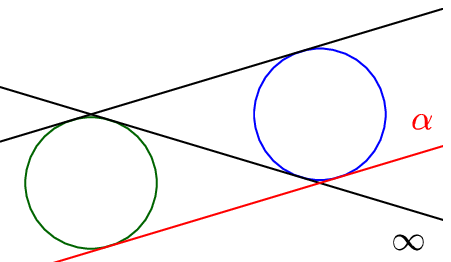}} \\ X тип, прямая
\end{minipage}
\vfill
\begin{minipage}[h]{0.3\linewidth}
\center{\includegraphics[width=1\linewidth]{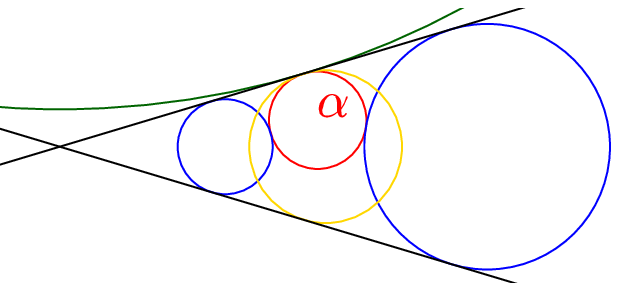}} \\ VIII тип
\end{minipage}
\hfill
\begin{minipage}[h]{0.3\linewidth}
\center{\includegraphics[width=1\linewidth]{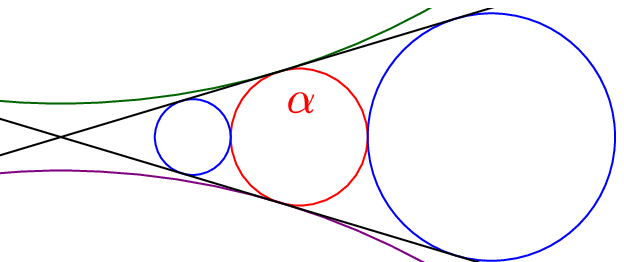}} \\ IX тип
\end{minipage}
\hfill
\begin{minipage}[h]{0.3\linewidth}
\center{\includegraphics[width=1\linewidth]{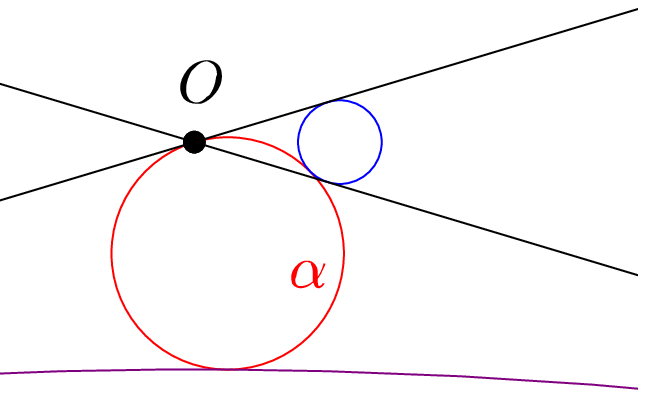}} \\ X тип, окружность
\end{minipage}
\caption{I-типы}
\label{ris:itypes}
\end{figure}

\begin{lemma}[об I-типах]\label{lem:itypes}
Даны две пересекающиеся прямые, а также окружность или прямая $\alpha$, отличная от них. Тогда, в зависимости от положения на плоскости, $\alpha$ принадлежит к одному из типов согласно Табл.~\ref{tab:itypes} (см. Рис.~\ref{ris:itypes}; мы их называем \emph{I-типами}).
\end{lemma}

\begin{table}[h!]
\begin{center}
\begin{tabular}{|c|c|c|c|c|}
\hline
Тип & Метка & Распределение решений для $\alpha$ & Доп. решения & Всего решений для $\alpha$ \\
\hline
I & $III_{1}$ & 2-2-2-2 & --- & 8 \\
\hline
II & $III_{2}$ & 4-2-0-2 & --- & 8 \\
\hline
III & $II$ & 2-2-0-0 & --- & 4 \\
\hline
IV & $I$ & 4-0-0-0 & --- & 4 \\
\hline
V & $[III]_{1}$ & 2-1-0-1 & $O$ или $\infty$ & 5 \\
\hline
VI & $[III]_{2}$ & 0-0-0-0 & $O$ и $\infty$ & 2 \\
\hline
VII & $IIT$ & 3-1-0-2 & --- & 6 \\
\hline
VIII & $IT$ & 3-1-0-0 & --- & 4 \\
\hline
IX & $ITT$ & 2-1-0-1 & --- & 4 \\
\hline
X & $[IIT]$ & 1-0-0-1 & $O$ или $\infty$ & 3 \\
\hline
\end{tabular}
\end{center}
\caption{I-типы}
\label{tab:itypes}
\end{table}

\textbf{Доказательство леммы~\ref{lem:itypes}} осуществляется одинаковым образом для всех типов посредством непрерывного увеличения допустимой окружности в одном из секторов до момента касания с $\alpha$. Более подробное рассуждение см. в~\cite{bib:general}.

\begin{remark}
Количество решений указано без учета кратности.
\end{remark}
\begin{remark}
Типы I и II имеют одинаковую метку, но разные распределения. В~I типе $\alpha$ разделяет $O$ и $\infty$, а в II --- нет. То же относится и к типам V и VI. В~V типе все три окружности имеют одну общую точку, а в VI --- две.
\end{remark}
\begin{remark}
Нетрудно заметить, что окружность II типа всегда пересекает обе стороны некоторого сектора. В~таком случае мы говорим, что эта окружность \emph{пересекает} этот сектор.
\end{remark}

\textbf{2. Зафиксированы две концентрические окружности.}

Под \emph{внутренней окружностью $\omega$} будем подразумевать окружность с меньшим радиусом, а под \emph{внешней $\Omega$} --- с большим. Все окружности, касающиеся $\omega$ и $\Omega$, разбиваются на два типа, которые мы условно обозначим A и B. Окружности \emph{типа A} касаются $\omega$ внешним образом, а \emph{типа B} --- внутренним. Их тоже логично назвать \emph{допустимыми}. Прямую, соединяющую центр $\omega$ и $\Omega$ с центром окружности $\alpha$, мы называем \emph{диаметром} (а луч --- \emph{радиусом}), на котором лежит $\alpha$. \emph{Противоположными} называются окружности, лежащие на одном диаметре, но по разные стороны от центра. Заметим, что любые две окружности типа B пересекаются.

Теперь снова добавим еще одну окружность $\alpha$, но теперь потребуем, чтобы она не имела общих точек с $\omega$ и $\Omega$. Найдем количество допустимых окружностей, касающихся $\alpha$.

\begin{figure}[h!]
\center{\includegraphics{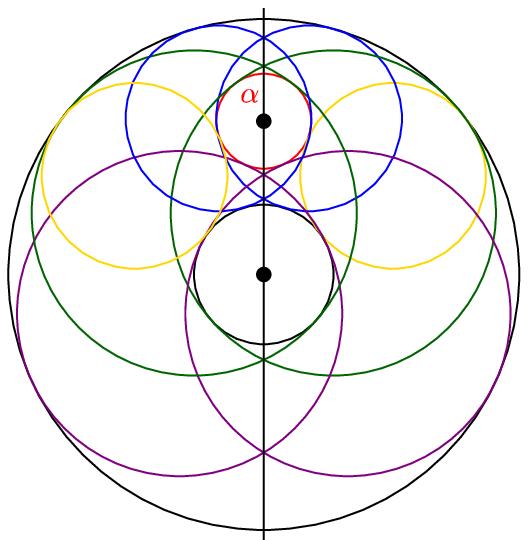}}
\caption{К лемме~\ref{lem:concentr}}
\label{ris:concentr}
\end{figure}

\begin{lemma}\label{lem:concentr}
Даны две концентрические окружности и окружность или прямая $\alpha$, не имеющая общих точек с данными. Тогда полученная тройка либо не имеет решений, либо имеет 8 решений (как на Рис.~\ref{ris:concentr}).
\end{lemma}

\textbf{Доказательство леммы~\ref{lem:concentr}.} Возможны следующие случаи:

\begin{enumerate}
  \item{$\alpha$ лежит внутри $\omega$. Тогда $\omega$ разделяет $\alpha$ и $\Omega$.}
  \item{$\alpha$ лежит снаружи $\Omega$. Тогда $\Omega$ разделяет $\alpha$ и $\omega$.}
  \item{$\alpha$ лежит внутри $\Omega$, но $\omega$ лежит внутри $\alpha$. Тогда $\alpha$ разделяет $\omega$ и $\Omega$.}
  \item{$\alpha$ лежит внутри $\Omega$, но $\omega$ не лежит внутри $\alpha$, и $\alpha$ не лежит внутри $\omega$.}
\end{enumerate}

Случаи 1) --- 3) соответствуют метке $S$, и тогда решений нет. Случай 4) соответствует пустой метке ($\varnothing$), и тогда есть 8 решений, причем 4 из них имеют тип A, и 4 --- тип B. Доказательство осуществляется посредством вращения потенциального решения до момента касания с $\alpha$.\qed

Обратим внимание, что полученные 8 решений разбиваются на пары симметричных относительно радиуса, на котором лежит $\alpha$. Очевидно также, что радиус $\alpha$ в этом случае меньше радиуса допустимых окружностей. Еще одно полезное наблюдение состоит в том, что только в одной паре симметричных решений окружности могут не пересекаться. Это пара решений типа A, которые касаются $\alpha$ внешним образом. Поэтому, в частности, среди любых 6 решений можно выбрать пару симметричных пересекающихся.

\textbf{3. Зафиксированы две параллельные прямые.}

Будем подразумевать направление фиксированных прямых \emph{горизонтальным}, а перпендикулярное --- \emph{вертикальным}. \emph{Допустимыми} являются окружности равного радиуса, касающиеся фиксированных прямых, а также все горизонтальные прямые. \emph{Распределением} здесь будем называть запись вида $x$+$y$, где $x$ обозначает количество допустимых окружностей, а $y$ --- количество допустимых прямых.

\begin{figure}
\begin{minipage}[h]{0.29\linewidth}
\center{\includegraphics[width=1\linewidth]{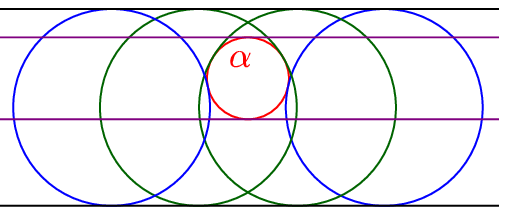}} \\ I тип
\end{minipage}
\hfill
\begin{minipage}[h]{0.29\linewidth}
\center{\includegraphics[width=1\linewidth]{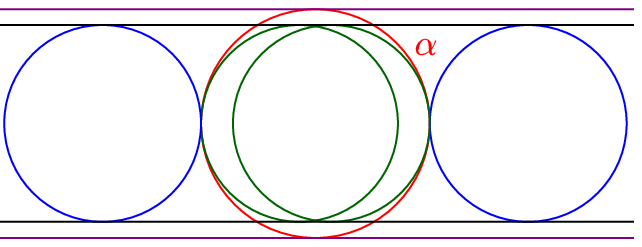}} \\ II тип
\end{minipage}
\hfill
\begin{minipage}[h]{0.29\linewidth}
\center{\includegraphics[width=1\linewidth]{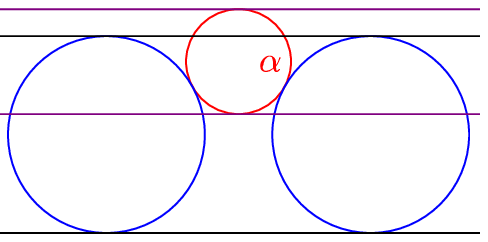}} \\ III тип
\end{minipage}
\vfill
\begin{minipage}[h]{0.29\linewidth}
\center{\includegraphics[width=1\linewidth]{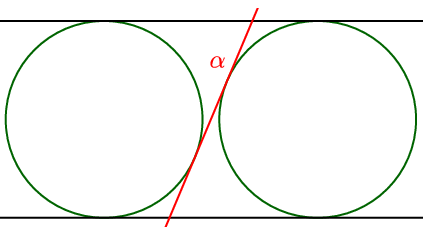}} \\ IV тип
\end{minipage}
\hfill
\begin{minipage}[h]{0.29\linewidth}
\center{\includegraphics[width=1\linewidth]{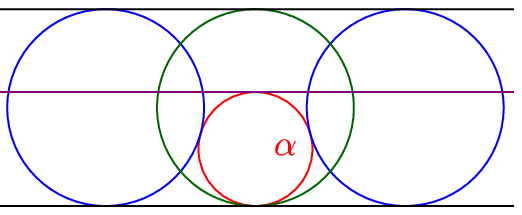}} \\ V тип
\end{minipage}
\hfill
\begin{minipage}[h]{0.29\linewidth}
\center{\includegraphics[width=1\linewidth]{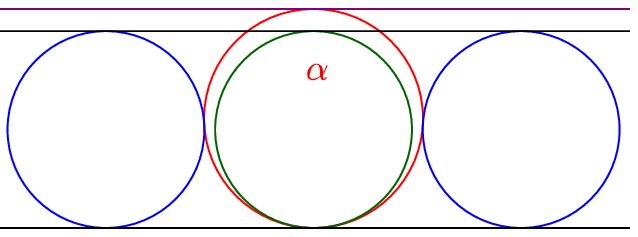}} \\ VI тип
\end{minipage}
\vfill
\begin{minipage}[h]{0.29\linewidth}
\center{\includegraphics[width=1\linewidth]{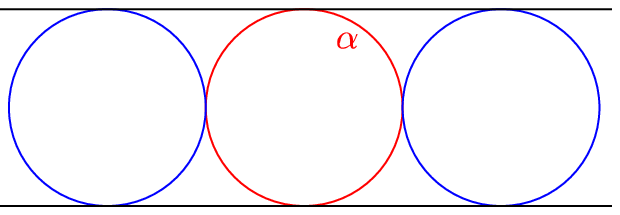}} \\ VII тип
\end{minipage}
\hfill
\begin{minipage}[h]{0.29\linewidth}
\center{\includegraphics[width=1\linewidth]{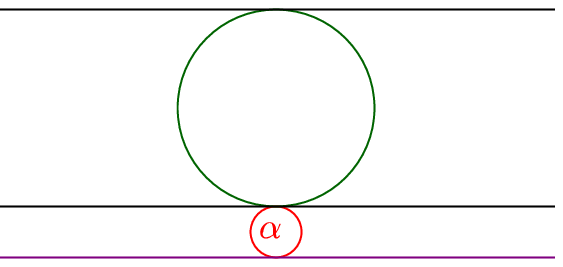}} \\ VIII тип
\end{minipage}
\hfill
\begin{minipage}[h]{0.29\linewidth}
\center{\includegraphics[width=1\linewidth]{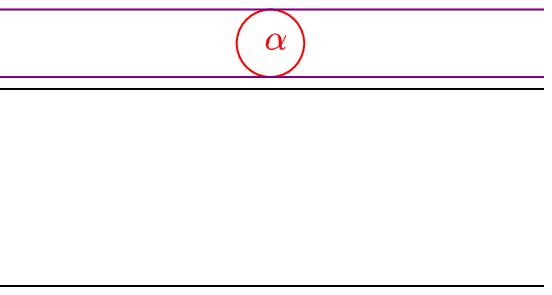}} \\ IX тип
\end{minipage}
\vfill
\center{\includegraphics[width=0.29\linewidth]{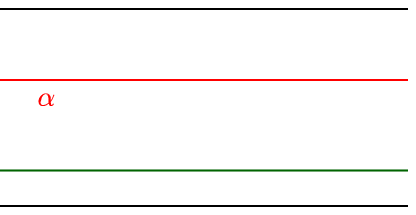}} \\ X тип
\caption{T-типы}
\label{ris:ttypes}
\end{figure}

\begin{lemma}[о T-типах]\label{lem:ttypes}
Даны две параллельные прямые и окружность или прямая $\alpha$. Тогда, в зависимости от положения на плоскости, $\alpha$ принадлежит к одному из типов согласно Табл.~\ref{tab:ttypes} (см. Рис.~\ref{ris:ttypes}; мы называем их \emph{T-типами}).
\end{lemma}

\begin{table}[h!]
\begin{center}
\begin{tabular}{|c|c|c|c|}
\hline
Тип & Метка & Распределение решений для $\alpha$ & Всего решений для $\alpha$ \\
\hline
I & $T$ & 4+2 & 6 \\
\hline
II & $IIT$ & 4+2 & 6 \\
\hline
III & $IT$ & 2+2 & 4 \\
\hline
IV & $[IIT]$ & 2+0 & 2 \\
\hline
V & $TT$ & 3+1 & 4 \\
\hline
VI & $ITT$ & 3+1 & 4 \\
\hline
VII & $TTT$ & 2+0 & 2 \\
\hline
VIII & $STT$ & 1+1 & 2 \\
\hline
IX & $ST$ & 0+2 & 2 \\
\hline
X & $[TTT]$ & 0+$\infty$ & $\infty$ \\
\hline
\end{tabular}
\end{center}
\caption{T-типы}
\label{tab:ttypes}
\end{table}

\textbf{Доказательство леммы~\ref{lem:ttypes}} осуществляется одинаковым образом для всех типов посредством непрерывного перемещения допустимой окружности между параллельными прямыми до момента касания с $\alpha$.

\begin{remark}
Как и в лемме~\ref{lem:itypes}, количество решений указано без учета кратности.
\end{remark}

\section{Парные решения}
Пусть на плоскости дано несколько окружностей и окружность или точка, касающаяся всех данных (решение). Сопоставим каждой данной окружности знак <<$+$>>, если она касается решения внешним образом и знак <<$-$>>, если внутренним (касание с точкой допускает любой знак). Полученное множество знаков мы называем \emph{комбинацией}. Комбинации называются \emph{противоположными}, если любая окружность в этих комбинациях имеет разный знак. Решение, соответствующее либо данной комбинации знаков, либо противоположной, назовем \emph{подходящим}. Два различных решения, подходящих для данной комбинации знаков назовем \emph{парными}.

\begin{lemma}\label{lem:algebraic}
Пусть дано множество окружностей $G$ и множество окружностей или точек $S$, для которых выполнены следующие условия:
\begin{enumerate}
  \item{В множестве $G$ хотя бы 3 окружности, $G$ --- невырожденное и каждой окружности из $G$ сопоставлен знак.}
  \item{Любая окружность или точка из $S$ является подходящей для данной комбинации знаков в $G$.}
\end{enumerate}
Тогда в $S$ не больше двух элементов.
\end{lemma}
Иными словами, если исходных окружностей хотя бы 3, то парных решений действительно не больше 2-х.

\textbf{Доказательство леммы~\ref{lem:algebraic}.} Мы приведем рассуждение, похожее на \cite[стр.~151-153]{bib:courant}, однако несколько более подробное и завершенное.

Выберем из $G$ три окружности и обозначим их множество через $G'$. Пусть $(x_i, y_i)$ --- координаты их центров, $r_i$ --- их радиусы, а $s_i=\pm1$ --- их знаки $(i = 1, 2, 3)$. Ясно, что если все окружности в $G'$ концентрические, то решений для тройки $G'$ нет. В~противном случае можно выбрать такую систему координат, в которой в каждом из наборов $(x_1, x_2, x_3)$ и $(y_1, y_2, y_3)$ есть хотя бы два различных числа (пусть для определенности $x_1\ne x_2$ и $x_1\ne x_3$). Рассмотрим систему с неизвестными $x, y, r$:

\begin{eqnarray}
(x-x_1)^2+(y-y_1)^2=(r+s_1r_1)^2\\
(x-x_2)^2+(y-y_2)^2=(r+s_2r_2)^2\\
(x-x_3)^2+(y-y_3)^2=(r+s_3r_3)^2
\end{eqnarray}

Заметим, что решения этой системы находятся во взаимно-однозначном соответствии с подходящими окружностями для $G'$. Действительно, рассмотрим любое подходящее решение с центром $(x_0, y_0)$ и радиусом $r_0$ (в случае точки $r_0=0$). Если оно соответствует знакам $s_i$, то $(x_0, y_0, r_0)$ является решением системы, а если оно соответствует знакам $-s_i$, то $(x_0, y_0, -r_0)$ является решением системы. Обратное соответствие также достаточно очевидно.

Отсюда ясно, что достаточно получить, что данная система имеет либо не более двух решений, либо бесконечно много решений. Дейcтвительно, если решений не больше двух, то и окружностей в $S$ не больше двух, что и требовалось. Если же решений бесконечно много, то и решений для $G'$ бесконечно много. Но бесконечным числом решений обладают только вырожденные тройки окружностей, а это противоречит условию 1).

Заменим уравнения $(2)$ и $(3)$ на разности $(2’)=(1)-(2)$ и $(3’)=(1)-(3)$, которые являются линейными уравнениями относительно $x, y, r$. В~силу нашего выбора системы координат, в них присутствуют как $x$, так и $y$, с ненулевыми коэффициентами хотя бы в одном из полученных уравнений.

Полученная система уравнений $(1), (2’), (3’)$ равносильна исходной. Рассмотрим $(2’), (3’)$ как систему линейных уравнений относительно $x$ и $y$. Заметим, что коэффициенты при $x$ и $y$ в этих уравнениях зависят только от координат центров выбранных окружностей из $G$, а их свободный член может зависеть от $r$ (однако только линейно). Возможны следующие случаи:
\begin{enumerate}
  \item{Коэффициенты при $x$ и $y$ не пропорциональны. Тогда $x$ и $y$ выражаются как линейные функции от $r$. Подставляя полученные выражения для $x$ и $y$ в $(1)$, получим уравнение на $r$, степени не более двух. Оно имеет либо не более двух, либо бесконечно много решений. Так как каждому значению $r$ однозначно соответствуют значения $x$ и $y$, лемма в этом случае доказана.}
  \item{Коэффициенты при $x$ и $y$ пропорциональны. Тогда, домножая одно из этих уравнений на число и вычитая его из другого, можно получить линейное уравнение на $r$. Разберем случаи в зависимости от количества его решений:}
  \begin{enumerate}
    \item{Это уравнение не имеет решений. Тогда решений не имеет и исходная система.}
    \item{Это уравнение имеет ровно одно решение. При соответствующем значении $r$ уравнения $(2')$ и $(3')$ совпадут. Значит, после подстановки полученного значения $r$ во все уравнения, исходная система равносильна системе из двух уравнений $(1)$ и $(2')$. Выразим из $(2')$ $y$ через $x$ и подставим его в $(1)$. Получится уравнение на $x$ степени не более двух. Оно имеет либо не более двух, либо бесконечно много решений. Так как каждому значению $x$ однозначно соответствует значение $y$, лемма доказана и в этом случае.}
    \item{Это уравнение имеет бесконечно много решений. Тогда уравнения $(2')$ и $(3')$ равносильны и из них можно оставить одно, из которого можно выразить $x$ как линейную функцию от $y$ и $r$. Подставив данное выражение для $x$ в $(1)$ получим уравнение не более чем второй степени от $y$ и $r$. На плоскости с осями $y, r$ оно может задавать пустое множество, одну точку, одну или две прямые, коническое сечение, а также всю плоскость. Таким образом, оно имеет либо не более одного, либо бесконечно много решений. Для получения решения исходной системы по $y, r$ достаточно их подставить в полученное линейное выражение для $x$. Значит, и вся система имеет либо не более одного, либо бесконечно много решений, что и требовалось.\qed}
  \end{enumerate}
\end{enumerate}
\begin{corollary}
Если центры окружностей из $G$ лежат на одной прямой, то парные решения симметричны относительно этой прямой.
\end{corollary}
\textbf{Доказательство следствия.} Предположим противное, т.~е. что нашлись не симметричные парные решения $\alpha$ и $\beta$. Тогда рассмотрим окружность $\alpha'$, симметричную относительно линии центров $G$. Ясно, что она также является подходящей и не совпадает с $\alpha$ и $\beta$. Получаем противоречие с леммой~\ref{lem:algebraic}\qed

\section{Доказательство теоремы~1}
\textbf{Доказательство}. Предположим противное условию теоремы, т.~е. что нашлись 4 объекта (не все из которых касаются в одной точке) с 7-ю или более решениями.

Предположим, что среди исходных объектов есть точка. Известно (см., например, \cite[стр.~100]{bib:general}), что если в тройке объектов есть точка, то количество решений не больше 4-х, за исключением случая вырожденной тройки. Поэтому в нашем случае любая тройка объектов, содержащая точку, вырождена, так как у нас есть 7 решений. Ясно, что тогда все объекты касаются в этой точке, что исключено условием.

Итак, все объекты --- окружности или прямые. Рассмотрим любую тройку исходных объектов. Согласно леммам~\ref{lem:itypes},~\ref{lem:concentr} и~\ref{lem:ttypes}, если для тройки объектов есть больше 6-и решений, то эта тройка имеет одну из меток $\varnothing$, III или [TTT]. Докажем, что последний случай невозможен. Действительно, достаточно рассмотреть любую другую тройку исходных объектов. Среди них точно какие-то два касаются, а значит эта тройка не может иметь метку $\varnothing$ или III. Поэтому она тоже имеет метку [TTT]. Тогда очевидно, что все объекты касаются в одной точке, что исключено формулировкой теоремы.

Таким образом, среди исходных объектов нет точек и пар касающихся.

Теперь рассмотрим два случая.
\paragraph{Первый случай: из исходных окружностей или прямых какие-то две пресекаются.}
Переведем инверсией любые две пересекающиеся окружности или окружность и прямую в пересекающиеся прямые. Для краткости здесь и далее образы исходных объектов при одной или нескольких инверсиях называются просто исходными объектами, а образы решений --- просто решениями.

Кроме двух пересекающихся прямых есть еще две исходные обобщенные окружности. Обозначим их $\alpha$ и $\beta$.

По лемме~\ref{lem:itypes}. для объектов не I или II типа уже по отдельности существует не больше 5-и решений, поэтому $\alpha$ и $\beta$ могут быть только объектами этих типов. Если они имеют разные типы, то их распределения равны 2-2-2-2 и 0-2-4-2 и они пересекаются не больше чем по 6-и окружностям т.~е. всего есть не более 6-и решений. Поэтому $\alpha$ и $\beta$ обязательно имеют одинаковые типы.

Пусть обе обобщенные окружности имеют тип I.

\begin{lemma}\label{lem:inter-1}
Если две окружности типа I имеют два общих решения в одном секторе, то они имеют общую точку в этом секторе.
\end{lemma}

\begin{figure}[h]
\begin{minipage}[h]{0.45\linewidth}
\center{\includegraphics[width=1\linewidth]{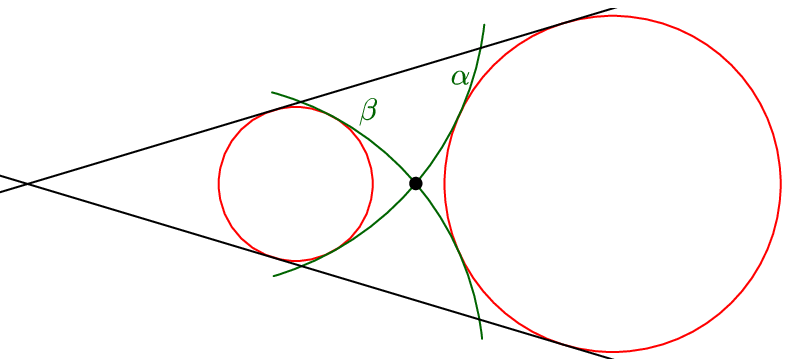}}
\end{minipage}
\hfill
\begin{minipage}[h]{0.45\linewidth}
\center{\includegraphics[width=1\linewidth]{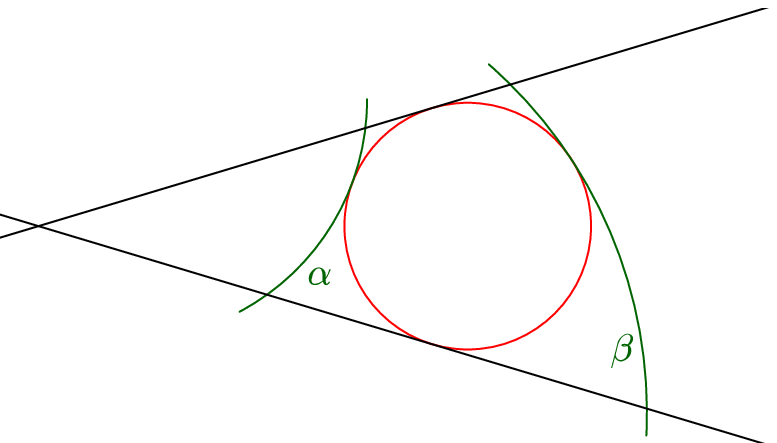}}
\end{minipage}
\caption{К лемме~\ref{lem:inter-1}}
\label{ris:inter-1}
\end{figure}

\textbf{Доказательство леммы~\ref{lem:inter-1}.} Заметим, что если $\alpha$ и $\beta$ не пересекаются в каком-то секторе, то они делят его на три области, из которых, возможно, лишь в одну можно вписать требуемую окружность-решение (Рис.~\ref{ris:inter-1} справа). Так как в нее, очевидно, можно вписать не более одной окружности, утверждение доказано.\qed

Так как у $\alpha$ и $\beta$ хотя бы 7 общих решений, то эти решения, очевидно, имеют распределение 2-2-2-1. Значит, согласно лемме~\ref{lem:inter-1} хотя бы в 3-х секторах $\alpha$ и $\beta$ должны пересекаться, что приводит к противоречию.

Пусть $\alpha$ и $\beta$ имеют тип II. Тогда, чтобы решений было хотя бы 7, необходимо, чтобы эти окружности пересекали один и тот же сектор, иначе уже по их распределениям (4-2-0-2 и 2-4-2-0, либо 4-2-0-2 и 0-2-4-2) можно сказать, что они имеют не более 4-х общих решений. С~этого сектора и начнем нумерацию. Тогда в I секторе не более 4-х решений, а значит во II и IV не меньше 3-х. Завершает доказательство первого случая следующая лемма.

\begin{lemma}\label{lem:symmetry}
Для двух окружностей $\alpha$ и $\beta$ II типа, пересекающих обе стороны сектора I, есть не более двух общих решений в секторах II и IV.
\end{lemma}
\textbf{Доказательство леммы~\ref{lem:symmetry}.} Предположим противное. Заметим, что $\alpha$ и $\beta$ касаются своих решений в секторах II и IV внешним образом. Центры этих решений лежат на одной прямой --- биссектрисе секторов II и IV. Ясно, что $\alpha$ и $\beta$ лежат по одну сторону от этой прямой, в частности, они не могут быть симметричны относительно нее. Противоречие с следствием из леммы~\ref{lem:algebraic}.\qed
\paragraph{Второй случай: никакие две окружности или прямые из исходных четырех не пересекаются.}
Переведем инверсией любые два из данных объектов в концентрические окружности. Образы двух оставшихся объектов снова назовем $\alpha$ и $\beta$.

По лемме~\ref{lem:concentr}, для $\alpha$ и $\beta$ отдельно существует 4 решения типа A и столько же типа B, причем каждое из этих множеств симметрично относительно радиуса, на котором расположена соответствующая окружность $\alpha$ или $\beta$. Так как у нас, согласно предположению, имеется не меньше 7-и решений, то хотя бы 4 из них одного типа. Далее нам потребуется несколько лемм.

\begin{figure}[h!]
\center{\includegraphics{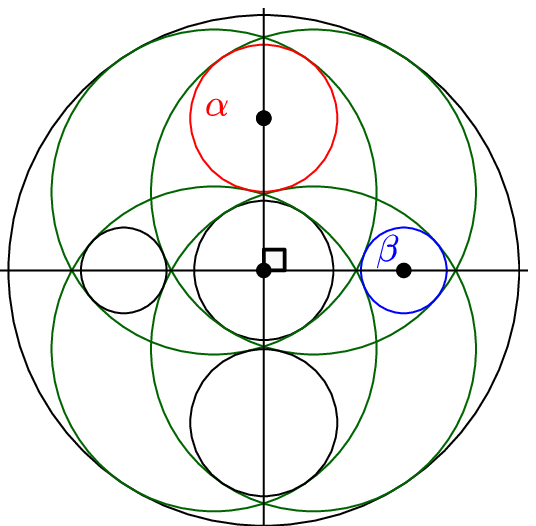}}
\caption{К лемме~\ref{lem:alternative}}
\label{ris:ortdiam}
\end{figure}

\begin{lemma}\label{lem:alternative}
Пусть окружности $\alpha$ и $\beta$ порождают 4 общих решения одного типа. Тогда либо $\alpha$ и $\beta$ лежат на одном диаметре, либо $\alpha$ и $\beta$ лежат на перпендикулярных диаметрах. В~последнем случае общие решения обязательно имеют тип B, а их центры образуют прямоугольник (см. Рис.~\ref{ris:ortdiam}).
\end{lemma}

\textbf{Доказательство леммы~\ref{lem:alternative}.} Пусть $\alpha$ и $\beta$ не лежат на одном диаметре (иначе доказывать нечего). Обозначим множество общих решений для $\alpha$ и $\beta$ через $X$. Поскольку окружности из $X$ являются решениями и для $\alpha$, и для $\beta$, то $X$ симметрично относительно радиусов $\alpha$ и $\beta$. Однако по нашему предположению, угол между диаметрами, на которых лежат $\alpha$ и $\beta$ ненулевой. Поэтому множество $X$ переходит в себя при некотором повороте на $\varphi>0^\circ$ относительно центра концентрических окружностей. Можно считать $\varphi\le180^\circ$.

Докажем, что $\varphi=90^\circ$ или $\varphi=180^\circ$. Рассмотрим произвольную окружность $\Gamma$ из $X$. При поворотах на 0, $\varphi$, 2$\varphi$, и т. д. она переходит в окружность из $X$. Но в $X$ ровно 4 окружности, поэтому $\Gamma$ может иметь от 2 до 4 образов при таких поворотах ($\Gamma$ не может перейти в себя, так как $0<\varphi\le180^\circ$). Эти случаи соответствуют поворотам на $180^\circ$, $120^\circ$ и $90^\circ$ соответственно. Однако если образов 3, то в $X$ остается еще одна окружность, не являющаяся образом $\Gamma$. У~нее тоже 3 образа, которые не совпадают с образами $\Gamma$. Получается, что в $X$ хотя бы 6 окружностей, а не 4 --- противоречие. Остаются варианты $180^\circ$ и $90^\circ$, которые реализуются, только если центры окружностей из $X$ образуют прямоугольник или квадрат соответственно.

Чтобы теперь доказать, что в $X$ не могут быть окружности типа A, нам понадобится следующая небольшая и достаточно очевидная лемма.

\begin{lemma}\label{lem:obvious}
Если для $\alpha$ и $\beta$ есть два общих противоположных решения $P$ и $Q$ типа A, то $\alpha$ и $\beta$ лежат на одном диаметре.
\end{lemma}

\begin{figure}[h!]
\center{\includegraphics{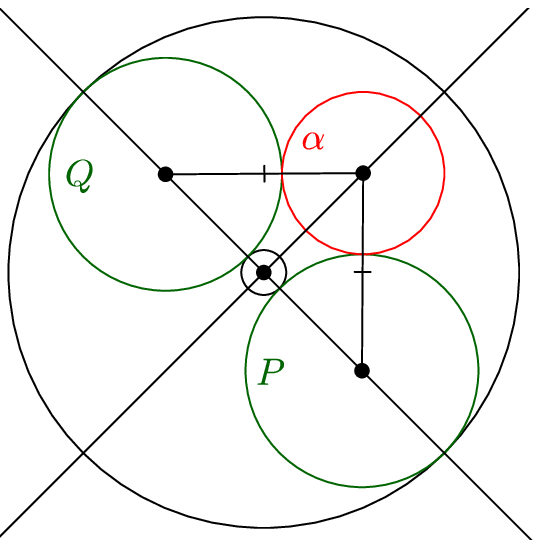}}
\caption{К лемме~\ref{lem:obvious}}
\label{ris:obvious}
\end{figure}
\textbf{Доказательство леммы~\ref{lem:obvious}.} Заметим, что $P$ и $Q$ не пересекаются, т.~к. они центрально-симметричны относительно точки, которой не содержат внутри себя. Так как радиус окружности $\alpha$ меньше радиуса окружности $P$ (в обычном смысле), то если $\alpha$ касается $P$ внутренним образом, то $\alpha$ лежит внутри $P$. Но т.~к. $P$ и $Q$ не пересекаются, то $\alpha$ тогда не может касаться $Q$. Значит, $\alpha$ и $\beta$ касаются $P$ и $Q$ внешним образом (см. Рис.~\ref{ris:obvious}). Но радиусы $P$ и $Q$ равны, поэтому ГМТ центров окружностей, касающихся внешним образом $P$ и $Q$, является серединным перпендикуляром отрезка между центрами $P$ и $Q$, который, очевидно, проходит через центр концентрических окружностей, что и требовалось.\qed

Продолжим доказательство леммы~\ref{lem:alternative}. Из леммы~\ref{lem:obvious} ясно, что в $X$ не могут быть окружности типа A, т.~к. окружности из $X$ делятся на две пары противоположных, поэтому $\alpha$ и $\beta$ должны лежать одновременно на двух диаметрах, что невозможно. Таким образом, в $X$ только окружности типа B. 

Так как центры окружностей из $X$ образуют прямоугольник, то кроме концентрических существуют еще 4 окружности, которые касаются всех окружностей из $X$ --- они вписаны в криволинейные симметричные четырехугольники, образованные окружностями из $X$ (см. Рис.~\ref{ris:ortdiam}). Больше окружностей, касающихся всех окружностей из $X$ быть не может по доказанному выше первому случаю основной теоремы. Действительно, в $X$ четыре попарно пересекающихся окружности, и их уже касаются 6 других (2 концентрические и 4 вписанные в криволинейные четырехугольники). Центры 4-х новых окружностей лежат на двух перпендикулярных диаметрах; $\alpha$ и $\beta$ должны быть двумя из этих четырех окружностей, т.~е. они должны лежать на двух перпендикулярных диаметрах.\qed

\begin{lemma}\label{lem:ortdiam}
Пусть окружности $\alpha$ и $\beta$ порождают 4 общих решения типа B и лежат на перпендикулярных диаметрах. Тогда они не могут порождать еще 3 общих решения типа A.
\end{lemma}
\textbf{Доказательство леммы~\ref{lem:ortdiam}.} Пусть это не так. Обозначим общие решения типа A через $P, Q, R$. Докажем, что среди них есть пара противоположных. Так как все они касаются $\alpha$, а всего среди решений для $\alpha$ четыре типа A, причем они разбиваются на две пары симметричных, то среди $P, Q, R$ две симметричны относительно радиуса, на котором лежит $\alpha$. Пусть это $P$ и $Q$. Рассмотрим теперь пару окружностей, аналогично симметричную относительно радиуса, на котором лежит $\beta$. Если это снова $P$ и $Q$, то они переходят друг в друга при двух симметриях, угол между осями которых составляет $90^\circ$, а тогда $P$ и $Q$ совпадают, что невозможно. Если же это $Q$ и $R$ (или, соответственно, $P$ и $R$), то тогда, так как угол между радиусами, на которых лежат $\alpha$ и $\beta$ равен $90^\circ$, то $P$ и $R$ (соответственно, $Q$ и $R$) противоположны --- отражение относительно двух перпендикулярных осей является центральной симметрией. Таким образом, среди $P, Q$ и $R$ есть пара противоположных. Следовательно, $\alpha$ и $\beta$ имеют два противоположных решения типа A, из чего с помощью леммы~\ref{lem:obvious} получаем, что они лежат на одном диаметре --- противоречие.\qed

Продолжим доказательство теоремы~1. Так как $\alpha$ и $\beta$ порождают 4 общих решения одного типа, то по лемме~\ref{lem:alternative} они лежат либо на одном диаметре, либо на перпендикулярных диаметрах. Но во втором случае $\alpha$ и $\beta$ порождают еще 3 общих решения типа A, что невозможно по лемме~\ref{lem:ortdiam}. Таким образом, $\alpha$ и $\beta$ лежат на одном диаметре.

Окружности $\alpha$ и $\beta$ порождают 7 общих решений, поэтому из них можно выбрать пару пересекающихся симметричных решений (см. наблюдение после леммы~\ref{lem:concentr}).

Сделаем инверсию в точке пересечения этих решений. Они перейдут в пересекающиеся прямые $a$ и $b$. Обратим внимание, что теперь эти прямые теперь будут играть роль исходных, а 4 исходные окружности образуют некоторое распределение. Следующая лемма завершает доказательство теоремы~1.

\begin{lemma}\label{lem:twosectors}
После инверсии исходные окружности образуют одно из распределений 4-0-0-0, 3-0-1-0 или 2-0-2-0 относительно прямых $a$ и $b$. Общее количество решений в этих случаях не превосходит 4, 2 и 6 соответственно. 
\end{lemma}
\textbf{Доказательство леммы~\ref{lem:twosectors}.} Центр инверсии лежит на радиусе, на котором лежат $\alpha$ и $\beta$ (так как точка пересечения двух симметричных окружностей лежит на оси симметрии). Поэтому центры всех 4-х исходных окружностей до инверсии лежали на одной прямой, проходящей через центр инверсии. Значит и после инверсии их центры лежат на одной прямой, а тогда все исходные окружности и должны лежать в двух противоположных секторах --- первая часть леммы, тем самым, доказана.

Обозначим 4 исходные окружности $\omega_1, \omega_2, \omega_3, \omega_4$ в том порядке, в котором они касаются любой из прямых $a$ или $b$.

Рассмотрим по отдельности все возможные распределения 4-х исходных окружностей. В~первом случае (распределение 4-0-0-0) по лемме~\ref{lem:itypes} оставшиеся решения должны быть II или IV типа. Но в IV типе среди 4-х решений (т.~е. исходных окружностей) заведомо есть пересекающиеся окружности --- те, что содержат окружность $\alpha$ внутри себя на Рис.~\ref{ris:itypes}. Значит, все оставшиеся решения имеют II тип. Они касаются $\omega_1$ и $\omega_4$ внешним образом, а $\omega_2$ и $\omega_3$ --- внутренним. По лемме~\ref{lem:algebraic}, где в качестве $G$ мы берем 4 исходные окружности, а в качестве $S$ --- решения II типа, имеем, что этих решений не более 2-х.

Для второго случая (3-0-1-0) вообще не существует подходящей окружности никакого типа, в чем нетрудно убедиться с помощью списка всех возможных распределений, приведенном в Табл.~\ref{tab:itypes}.

В третьем же случае (2-0-2-0) оставшиеся решения могут быть I или II типа.

\begin{figure}[h]
\begin{minipage}[h]{0.35\linewidth}
\center{\includegraphics[width=1\linewidth]{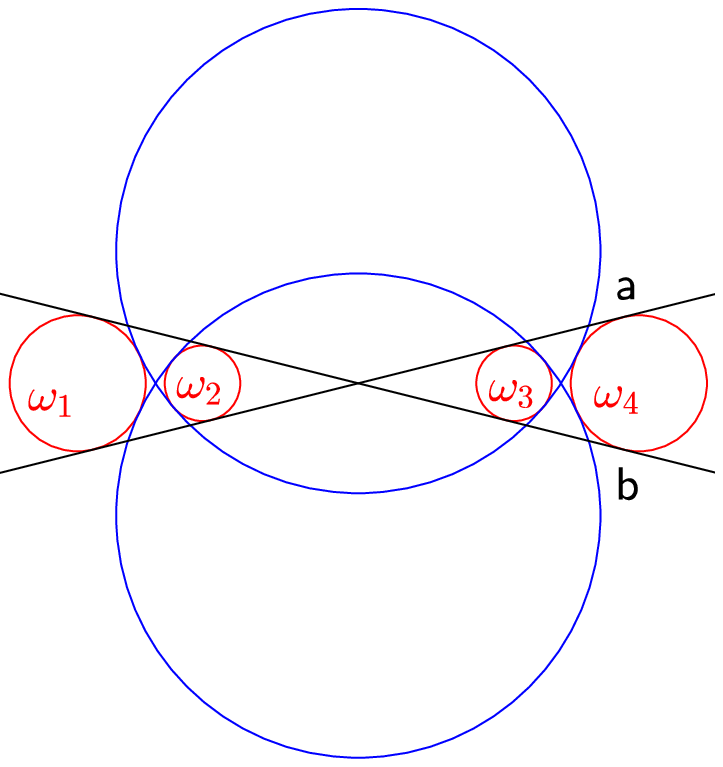} \\ Решения I типа}
\end{minipage}
\hfill
\begin{minipage}[h]{0.35\linewidth}
\center{\includegraphics[width=1\linewidth]{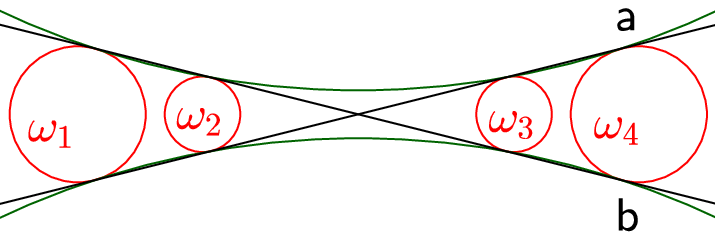} \\ Решения II типа}
\end{minipage}
\caption{Возможные пары решений при распределении 2-0-2-0}
\label{ris:final}
\end{figure}

Решения II типа должны пересекать сектора II и IV и касаться исходной четверки внешним образом. Применим для них лемму~\ref{lem:algebraic}, взяв в качестве $G$ четверку исходных окружностей, а в качестве $S$ --- все решения типа II. Получаем, что есть не более двух решений II типа (Рис.~\ref{ris:final} справа).

Решения I типа касаются $\omega_2$ и $\omega_3$ внутренним образом, а $\omega_1$ и $\omega_4$ --- внешним. Применим еще раз лемму~\ref{lem:algebraic}, взяв в качестве $G$ четверку исходных окружностей, а в качестве $S$ --- все решения типа I. Получаем, что есть не более двух решений I типа (Рис.~\ref{ris:final} слева).

Таким образом, общее количество решений типов I и II в этом случае не больше 4-х, т.~е. общее количество решений не больше 6-и и доказательство леммы~\ref{lem:twosectors}, а с ней и всей теоремы~1 завершено.\qed

\section{Доказательство теоремы~2}
\textbf{Доказательство}. Как и при доказательстве теоремы~1 предположим противное, т.~е. что нашлась пятерка объектов (не все из которых касаются в одной точке) с 5-ю или более решениями. Ясно, что среди исходных объектов нет точек (см. начало док-ва теоремы~1). Поскольку количество решений равно количеству исходных объектов, то решения можно рассматривать в качестве исходных объектов (и наоборот), поэтому среди решений тоже нет точек. Также пока предположим, что среди исходных объектов и решений нет касающихся.

Заметим, что среди имеющихся 10-и обобщенных окружностей есть пересекающиеся. Действительно, в противном случае переведем любые две непересекающиеся обобщенные окружности в концентрические. Тогда среди 5-и решений для любой третьей окружности достаточно найти два пересекающихся. Если есть два решения типа B, то они являются искомыми. В~противном случае есть все 4 решения типа A, среди которых есть пересекающиеся парные решения (см. наблюдение после леммы~\ref{lem:concentr}).

Рассмотрим соответствующую пятерку объектов в качестве исходной и переведем пересекающиеся объекты в прямые. Оставшиеся исходные объекты обозначим $\alpha$, $\beta$ и $\gamma$. Из леммы~\ref{lem:itypes} получаем, что каждый из них должен иметь I, II или V тип. Но если какой-то объект имеет V тип, то среди решений должна быть точка, что невозможно. Значит, каждый из объектов $\alpha$, $\beta$ и $\gamma$ имеет I или II тип. Соответственно, возможно 4 случая. Прежде чем разбирать их заметим, что две окружности II типа могут иметь больше 4-х общих решений только тогда, когда пересекают один и тот же сектор. Поэтому можно считать, что все встречающиеся ниже окружности II типа пересекают I сектор.

\boldmath\textbf{1. Все окружности $\alpha$, $\beta$ и $\gamma$ I типа.}\unboldmath\ Нам потребуются две новые леммы.
\begin{lemma}\label{lem:inter-2}
Пусть для окружностей $\alpha$, $\beta$ и $\gamma$ I типа есть общее решение в некотором секторе. Тогда среди них какие-то две имеют общую точку в этом секторе.
\end{lemma}
\begin{figure}[h]
\center{\includegraphics{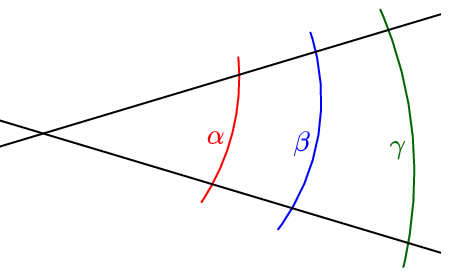}}
\caption{К лемме~\ref{lem:inter-2}}
\label{ris:inter-2}
\end{figure}
\textbf{Доказательство леммы~\ref{lem:inter-2}.} Предположим противное. Тогда дуги окружностей $\alpha$, $\beta$ и $\gamma$ делят сектор на четыре части, каждая из которых граничит не более чем c двумя из этих окружностей (см. Рис.~\ref{ris:inter-2}). Потенциальное решение должно быть вписано в одну из этих частей, что, в таком случае, невозможно.\qed

\begin{lemma}\label{lem:inter-3}
Пусть для окружностей $\alpha$, $\beta$ и $\gamma$ I типа есть два общих решения в некотором секторе. Тогда среди них какие-то две имеют две общие точки этом секторе.
\end{lemma}
\begin{figure}[h]
\begin{minipage}[h]{0.4\linewidth}
\center{\includegraphics[width=1\linewidth]{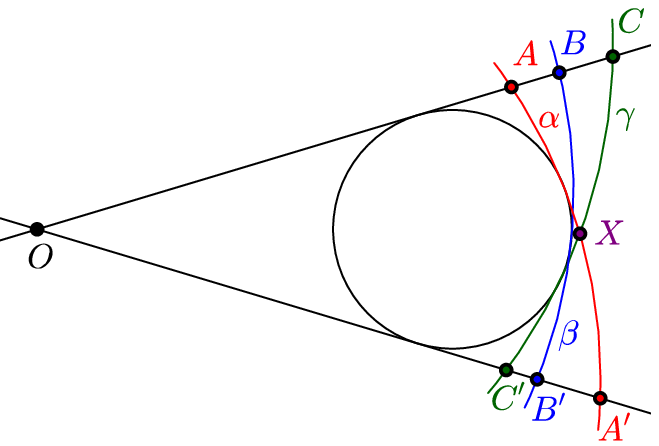} \\ $\beta$ разделяет $O$ и $X$}
\end{minipage}
\hfill
\begin{minipage}[h]{0.4\linewidth}
\center{\includegraphics[width=1\linewidth]{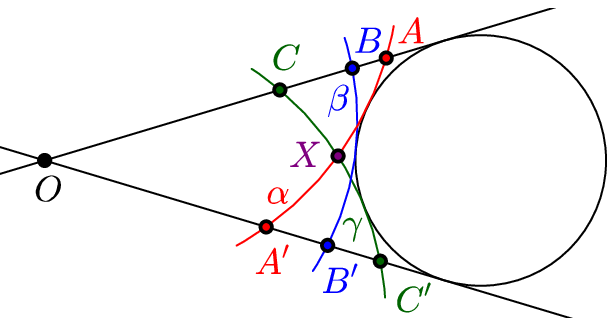} \\ $\beta$ не разделяет $O$ и $X$}
\end{minipage}
\caption{К лемме~\ref{lem:inter-3}}
\label{ris:inter-3}
\end{figure}
\textbf{Доказательство леммы~\ref{lem:inter-3}.} Предположим противное. Пусть дуги окружностей $\alpha$, $\beta$ и $\gamma$ пересекают одну из сторон данного сектора в точках $A$, $B$ и $C$, а вторую --- в точках $A'$, $B'$ и $C'$ соответственно. Пусть для определенности $B$ лежит на отрезке $AC$.

По лемме~\ref{lem:inter-1} дуги окружностей $\alpha$, $\beta$ и $\gamma$ имеют общую точку в рассматриваемом секторе причем, согласно нашему предположению, ровно одну. Но тогда порядок троек точек $O$, $A$, $B$ и $O$, $A'$, $B'$ должен отличаться, т.~е. если $A$ лежит на отрезке $OB$, то $B'$ лежит на отрезке $OA'$, и наоборот: если $B$ лежит на отрезке $OA$, то $A'$ лежит на отрезке $OB'$. То же верно и для двух других пар дуг окружностей ($\alpha$ и $\gamma$, $\beta$ и $\gamma$). Отсюда очевидно, что $B'$ лежит на $A'C'$.

Пусть теперь $X$ --- точка пересечения дуг $\alpha$ и $\gamma$. Возможно два случая, в зависимости от того, разделяет ли $\beta$ точки $O$ и $X$ (см. Рис.~\ref{ris:inter-3}). Нетрудно убедиться, что в обоих случаях $\alpha$, $\beta$ и $\gamma$ делят сектор на 7 частей, из которых только одна граничит со сторонами сектора и всеми тремя данными окружностями --- только в эту часть может быть вписана искомая окружность. Получаем противоречие, поскольку всего решений два.\qed

Продолжим разбор первого случая теоремы~2. Если окружности $\alpha$, $\beta$ и $\gamma$ имеют 5 общих решений, то их распределение изоморфно 2-2-1-1, 2-1-2-1 или 1-1-1-2. Пусть в распределении $k$ двоек и $4-k$ единиц, где $k=1,2$. Согласно лемме~\ref{lem:inter-2} в секторах с одним решением есть хотя бы 2 точки пересечения $\alpha$, $\beta$ и $\gamma$, а в секторах с двумя решениями их хотя бы 4 по леммам~\ref{lem:inter-1} и~\ref{lem:inter-3}. Поэтому общее число точек пересечения $\alpha$, $\beta$ и $\gamma$ не меньше чем $4k+(4-k)=3k+4\ge7$, а их на самом деле не более шести --- противоречие.

\boldmath\textbf{2. Все окружности $\alpha$, $\beta$ и $\gamma$ II типа.}\unboldmath\ Рассмотрим случаи в зависимости от количества решений в I секторе. Обозначим эти решения $\omega_i, i=1,2,3,4$ в том порядке, в котором они касаются любой из прямых.

\textit{2.1. В I секторе 4 решения.} Тогда $\alpha$, $\beta$ и $\gamma$ касаются $\omega_1$, $\omega_4$ внешним образом, а $\omega_2$, $\omega_3$ --- внутренним. Из леммы~\ref{lem:algebraic} для $G=\{\omega_1,\omega_2,\omega_3,\omega_4\}$ и $S=\{\alpha,\beta,\gamma\}$ получаем противоречие.

\textit{2.2. В I секторе 3 решения.} Рассмотрим окружность $\alpha$. Окружности $\omega_1,\omega_2,\omega_3$ являются тремя из 4-х решений для $\alpha$ в I секторе (см. Рис.~\ref{ris:itypes}). Простой перебор показывает, что если взять любые 3 решения из этих 4-х, то среднее из них касается $\alpha$ внутренним образом. Поэтому $\alpha$, $\beta$ и $\gamma$ касаются $\omega_2$ внутренним образом. Однако есть еще 2 решения в секторах II и IV, которые касаются $\alpha$, $\beta$ и $\gamma$ внешним образом. Применим лемму~\ref{lem:algebraic} взяв в качестве $G$ эти 2 решения и $\omega_2$, а в качестве $S$ --- $\alpha$, $\beta$ и $\gamma$, получая противоречие.

\textit{2.3. В I секторе не более 2-х решений.} Тогда в секторах II и IV хотя бы 3 решения, что противоречит лемме~\ref{lem:symmetry}.

\boldmath
\textbf{3. Две окружности II типа (например, $\alpha$ и $\gamma$) и одна окружность I типа ($\beta$).} Тогда за счет окружности I типа в I секторе не более 2-х решений, но тогда в II и IV секторах хотя бы 3 решения, что противоречит лемме~\ref{lem:symmetry}.

\textbf{4. Две окружности I типа ($\alpha$ и $\gamma$) и одна окружность II типа ($\beta$).} В этом случае возможны распределения 1-2-0-2 и 2-2-0-1. Разберем их отдельно.
\unboldmath

\textit{4.1. Распределение 1-2-0-2.}
\begin{lemma}\label{lem:inter-4}
Даны окружность $\alpha$ I типа и окружность $\beta$ II типа, пересекающая II сектор, причем $\alpha$ и $\beta$ имеют в I секторе 2 общих решения. Тогда они имеют общую точку в этом секторе.
\end{lemma}
\begin{figure}[h]
\center{\includegraphics{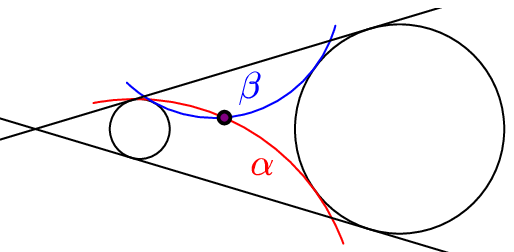}}
\caption{К лемме~\ref{lem:inter-4}}
\label{ris:inter-4}
\end{figure}
\textbf{Доказательство леммы~\ref{lem:inter-4}.} Заметим, что окружность I типа всегда разделяет свои решения в данном секторе. Поэтому и дуга $\alpha$ в I секторе разделяет свои решения. Рассмотрим дугу~$\beta$, лежащую во I секторе. Она имеет общие точки с обоими решениями (поскольку $\beta$ их касается), а значит линия между этими точками должна пересекать дугу $\alpha$ в этом секторе (в силу непрерывности окружности), что и требовалось.\qed

\begin{figure}[h!]
\center{\includegraphics{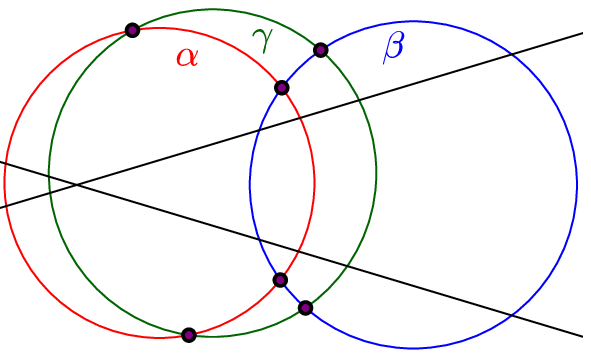}}
\caption{К доказательству случая 4.1}
\label{ris:5pieces}
\end{figure}

Продолжим доказательство случая~4.1 теоремы~2. Так как в секторах II и IV по 2 решения, то из леммы~\ref{lem:inter-4} получаем (после перенумерации секторов), что все точки пересечения $\beta$ с $\alpha$ и $\gamma$ находятся в этих секторах, а по лемме~\ref{lem:inter-1} окружности $\alpha$ и $\gamma$ также должны иметь общие точки в этих секторах. Таким образом, $\alpha$, $\beta$ и $\gamma$ попарно не имеют общих точек в секторе I. В~таком случае, они разбивают его на 5 частей (см. Рис.~\ref{ris:5pieces}), в каждую из которых можно вписать окружность, касающуюся лишь не более двух из $\alpha$, $\beta$ и $\gamma$ --- противоречие.

\begin{figure}[h!]
\center{\includegraphics[scale=0.9]{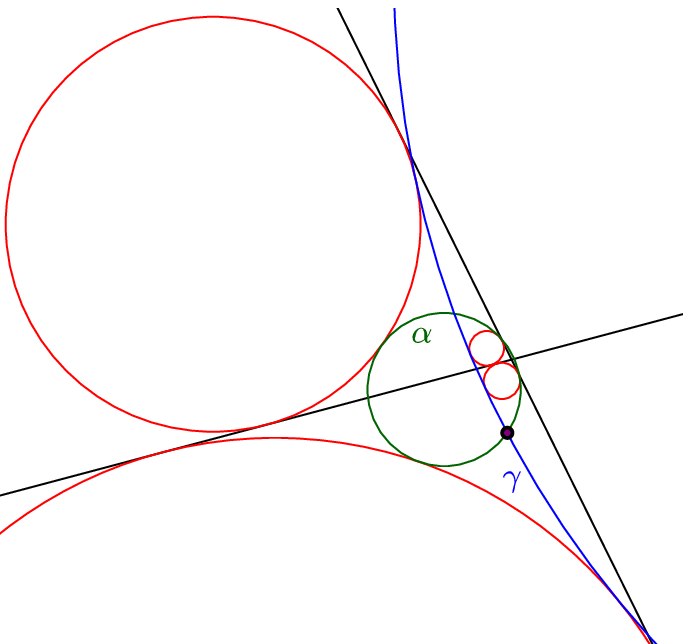}}
\caption{Инверсия в точке пересечения $\alpha$ и $\gamma$}
\label{ris:lastcase}
\end{figure}

\textit{4.2. Распределение 2-2-0-1.} Сделаем инверсию в точке пересечения $\alpha$ и $\gamma$. Тогда уже они перейдут в пересекающиеся прямые. Заметим, что $\alpha$ и $\gamma$ касались 2-х решений в секторах I и II внешним образом, а двух других --- внутренним, т.~е. два из этих решений лежали внутри $\alpha$ и $\gamma$, а два других --- вне (см. Рис.~\ref{ris:lastcase}). Но 4 области, на которые две пересекающиеся окружности делят плоскость, при инверсии в точке их пересечения перейдут в сектора соответствующих прямых. Поэтому получившееся распределение имеет в некоторых двух противоположных секторах не менее двух решений в каждом, поэтому распределение решений не может быть 2-2-0-1. Значит, этот случай нами уже рассмотрен и, следовательно, невозможен.

Остается рассмотреть случаи, когда среди исходных обобщенных окружностей или решений есть касающиеся.

Сделаем инверсию в точке касания этих обобщенных окружностей. Они перейдут в (горизонтальные) параллельные прямые. Заметим, что среди оставшихся трех исходных объектов (мы их обозначим, как обычно, $\alpha$, $\beta$ и $\gamma$) не может быть горизонтальной прямой. Действительно, тогда конфигурация исходных окружностей вырождена и все решения являются горизонтальными прямыми. Но для любого из остальных T-типов есть не более 2-х прямых-решений. Поэтому все исходные окружности являются параллельными прямыми, что противоречит условию теоремы.

Среди T-типов только два имеют 5 или более решений. Это типы I и II, их распределение одинаково и равно 4+2. Значит распределение 5-и решений может иметь вид 3+2 или 4+1. Разберем эти случаи.

\newpage\textbf{1. Распределение 4+1.} 
\begin{lemma}\label{lem:paralell}
Если окружности $\alpha$ и $\beta$ порождают 4 общих окружности-решения, то их центры лежат на одной вертикальной прямой.
\end{lemma}
\textbf{Доказательство леммы~\ref{lem:paralell}.} Предположим противное и обозначим множество из 4-х общих окружностей-решений через $X$. Это множество симметрично относительно вертикальных прямых, на которых лежат центры $\alpha$ и $\beta$. По нашему предположению, эти прямые не совпадают, поэтому $X$ переходит в себя при симметрии относительно каждой из некоторых двух параллельных прямых. Но тогда $X$ переходит в себя и при некотором параллельном переносе т.~е. является бесконечным --- противоречие.\qed

\begin{figure}[h]
\center{\includegraphics{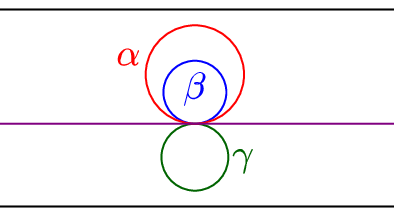}}
\caption{Распределение 4+1}
\label{ris:deg1}
\end{figure}

По этой лемме имеем, что центры $\alpha$, $\beta$ и $\gamma$ лежат на одной вертикальной прямой, а тогда они касаются прямой-решения в одной и той же точке, т. е. образуют вырожденную тройку (см. Рис.~\ref{ris:deg1}), что исключено выше.

\textbf{2. Распределение 3+2.} Заметим, что поскольку среди решений есть две параллельные прямые, то радиусы $\alpha$, $\beta$ и $\gamma$ равны половине расстояния между этими прямыми, а их центры лежат на одной горизонтальной прямой (см. Рис.~\ref{ris:deg2} слева). То же можно сказать и про три оставшихся решения-окружности (поскольку среди исходных объектов есть две параллельные прямые). Получаем два множества окружностей $X$ и $Y$, которые удовлетворяют следующим условиям:
\begin{figure}
\begin{minipage}[h]{0.4\linewidth}
\center{\includegraphics{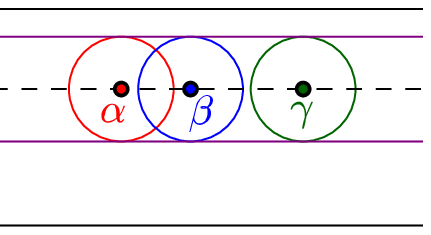} \\ Распределение 3+2}
\end{minipage}
\hfill
\begin{minipage}[h]{0.4\linewidth}
\center{\includegraphics{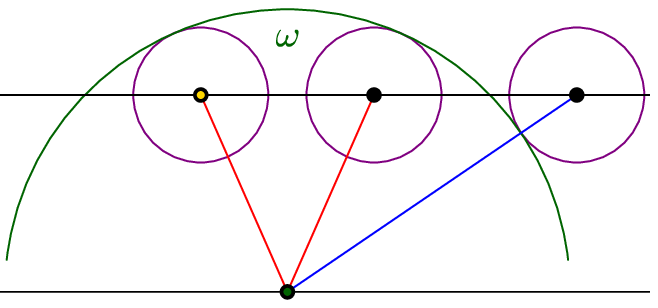} \\ Финальное рассуждение}
\end{minipage}
\caption{}
\label{ris:deg2}
\end{figure}
\begin{enumerate}
\item{В $X$ и $Y$ по три окружности.}
\item{В каждом из множеств окружности имеют равный радиус.}
\item{В каждом из множеств центры окружностей лежат на одной прямой, причем для $X$ и для $Y$ эти прямые параллельны (горизонтальны).}
\item{Любые две окружности из разных множеств касаются.}
\end{enumerate}
Мы докажем, что одновременное соблюдение этих условий невозможно.

Возьмем самую \emph{левую} окружность $\omega$ из $X\cup Y$. Расстояние между центром окружности из $X$ и центром окружности из $Y$ может принимать одно из двух значений, в зависимости от того, внутреннее касание или внешнее. Поэтому из трех таких расстояний для окружности $\omega$ найдутся два равных. Центры соответствующих окружностей образуют равнобедренный треугольник, одна из вершин которого окажется \emph{левее} $\omega$ (Рис.~\ref{ris:deg2} справа). Полученное противоречие завершает доказательство данного случая и всей теоремы~2.\qed

\section{Описание всех четверок окружностей с 6-ю решениями}
Итак, по теореме~1 для 4-х окружностей почти всегда существует не более 6-и решений. Оказывается, все четверки, для которых это количество достигается, допускают общее описание.

\begin{figure}[h!]
\begin{minipage}{0.24\linewidth}
\center{\includegraphics[width=1\linewidth]{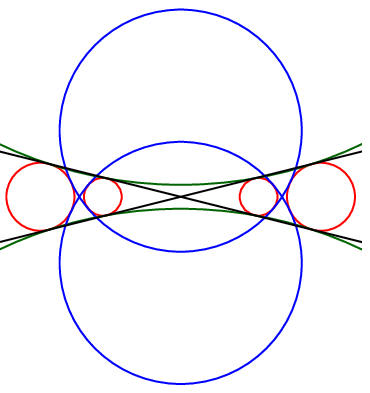}} \\ №1
\end{minipage}
\hfill
\begin{minipage}{0.24\linewidth}
\center{\includegraphics[width=1\linewidth]{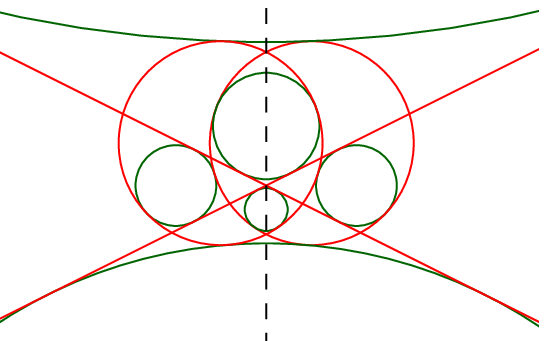}} \\ №2
\end{minipage}
\hfill
\begin{minipage}{0.24\linewidth}
\center{\includegraphics[width=1\linewidth]{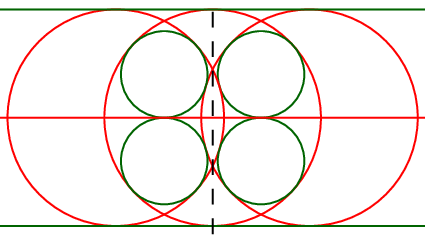}} \\ №3
\end{minipage}
\hfill
\begin{minipage}{0.24\linewidth}
\center{\includegraphics[width=1\linewidth]{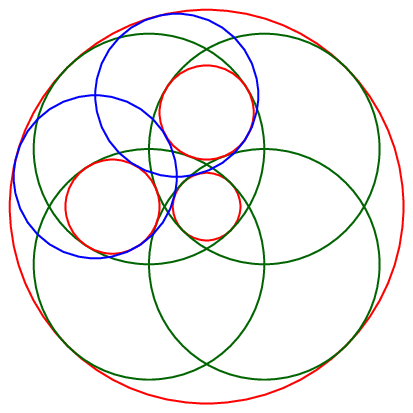}} \\ №4
\end{minipage}
\caption{Все четверки окружностей с 6-ю решениями}
\label{ris:forms}
\end{figure}

\begin{theorem}
Даны 4 различные обобщенные окружности, которые не все касаются в одной точке. Известно, что существует 6 окружностей, касающихся каждой из данных. Тогда данные 10 окружностей можно инверсиями и подобиями привести к одному из типов на Рис.~\ref{ris:forms}.
\end{theorem}

\begin{remark}[построение №1]
Рассмотрим пересекающиеся прямые. Далее возьмем 2 допустимые окружности в I секторе и 2 симметричные им относительно $O$ (это исходные окружности). Остается достроить решения I и II типов для полученной четверки (они, очевидно, существуют).
\end{remark}
\begin{remark}[построение №2]
Рассмотрим пересекающиеся прямые (это 2 исходные окружности). Далее возьмем допустимую окружность $\omega$ в I секторе и симметричную ей относительно $O$ окружность $\omega'$ в III секторе. Затем возьмем окружность I типа, касающуюся $\omega$ внешним образом, а $\omega'$ --- внутренним, и симметричную ей относительно биссектрисы секторов II и IV (это еще 2 исходные окружности). Остается достроить 4 решения в секторах II и IV, которые, очевидно, существуют.
\end{remark}
\begin{remark}[построение №3]
Рассмотрим параллельные прямые на расстоянии 2 и параллельную им прямую на расстоянии $1$ от них (это 1-я исходная обобщенная окружность). Далее возьмем 3 допустимые окружности, причем расстояние между центрами соседних из них равно $2/\sqrt{5}$ (это оставшиеся 3 исходные окружности). Достроим остальные 4 решения.
\end{remark}
\begin{remark}[построение №4]
Рассмотрим концентрические окружности с отношением радиусов $3\pm2\sqrt{2}$ (это 2 исходные окружности). Далее возьмем 4 окружности типа B с центрами в вершинах квадрата. Оставшиеся 2 исходные окружности касаются каждой из построенных окружностей типа B и лежат на перпендикулярных диаметрах. Остается достроить 2 решения типа A.
\end{remark}

\textbf{Доказательство теоремы~3.} Снова ясно, что среди исходных обобщенных окружностей нет точек (см. начало доказательства теоремы~1). Докажем, что среди них нет и касающихся. Пусть это не так, тогда переведем касающиеся объекты в параллельные прямые. Останутся два объекта $\alpha$ и $\beta$, которые должны иметь I или II тип. Все порожденные ими решения общие. Среди них есть параллельные прямые, поэтому линия центров $\alpha$ и $\beta$ горизонтальна. Однако среди решений есть и 4 окружности, поэтому по лемме~\ref{lem:paralell} линия центров $\alpha$ и $\beta$ вертикальна --- противоречие.

Как и в теореме~1, рассмотрим случаи пересекающихся и непересекающихся исходных обобщенных окружностей. 
\paragraph{Первый случай: из исходных окружностей или прямых какие-то две пресекаются.}
Переведем их инверсией в пересекающиеся прямые. Заметим, что среди I-типов, которые не содержат касающихся объектов, только типы I и II имеют не менее 6 решений. Поэтому оставшиеся окружности $\alpha$ и $\beta$ имеют I или II тип, причем теперь все 3 случая возможны. Разберем их по отдельности.

\begin{figure}[h!]
\center{\includegraphics{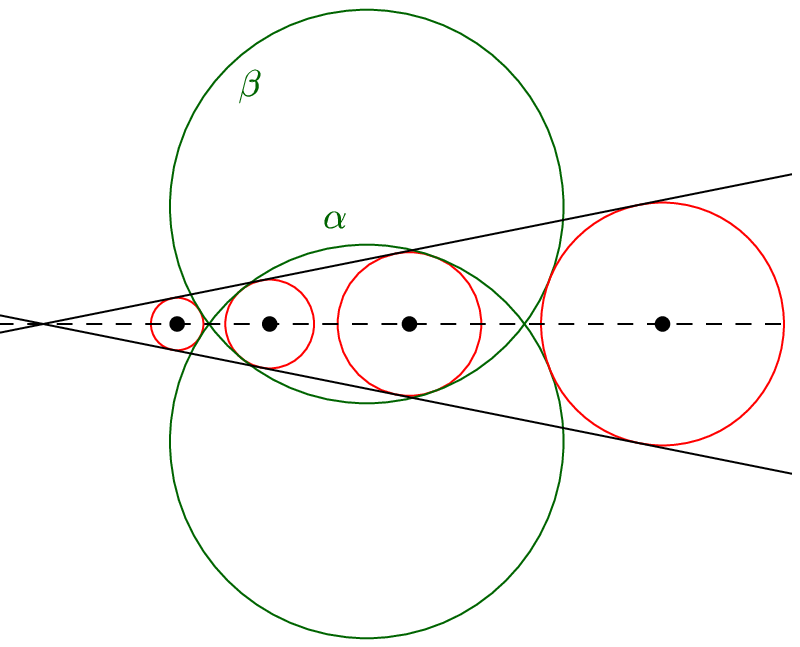}}
\caption{Обе исходные окружности II типа}
\label{ris:bothII}
\end{figure}

\boldmath\textbf{1. Обе окружности $\alpha$ и $\beta$ имеют II тип.}\unboldmath\ Тогда $\alpha$ и $\beta$ пересекают один и тот же сектор --- пусть это I сектор. Тогда по лемме~\ref{lem:symmetry} в секторах II и IV не более 2-х решений, а значит в I секторе 4 решения, причем $\alpha$ и $\beta$ касаются их в одной и той же комбинации знаков <<$+--+$>>, считая от $O$ (см. Рис.~\ref{ris:itypes}). Поэтому по следствию из леммы~\ref{lem:algebraic} они симметричны относительно биссектрисы секторов I и III (см. Рис.~\ref{ris:bothII}).

Очевидно, что центр любой окружности, касающейся $\alpha$ и $\beta$ внешним образом должен лежать на их оси симметрии (т.~к. радиусы $\alpha$ и $\beta$ равны) т.~е. на биссектрисе секторов I и III. Но центры решений в секторах II и IV, которые как раз касаются $\alpha$ и $\beta$ внешним образом, лежат на биссектрисе секторов II и IV --- противоречие. Значит, этот случай вообще невозможен.

\boldmath\textbf{2. Одна окружность I типа ($\alpha$) и одна окружность II типа ($\beta$).}\unboldmath\ Пусть $\beta$ пересекает I сектор. Тогда распределение решений равно 2-2-0-2. Сперва нам потребуется следующая (интересная и сама по себе) лемма.
\begin{lemma}\label{lem:beaty}
Окружность I типа $\alpha$ и окружность II типа $\beta$ порождают 4 общих решения с распределением 2-0-2-0. Тогда множество этих решений симметрично относительно точки пересечения прямых $O$.
\end{lemma}
\textbf{Доказательство леммы~\ref{lem:beaty}.} Обозначим решения через $\omega_i$ (в том порядке, в котором они касаются любой из прямых), а через $a_i$ обозначим длину внешней касательной из $O$ к $\omega_i$ ${(i=1,2,3,4)}$.

Рассмотрим инверсию с центром в $O$, оставляющую $\beta$ на месте (такая инверсия существует т.~к. $O$ лежит вне любой окружности II типа). При такой инверсии прямые также остаются на месте, поэтому решения должны перейти в решения. Ясно, что после инверсии любое решение остается в своем секторе. Следовательно, $\omega_1$ меняется с $\omega_2$, а $\omega_3$ меняется с $\omega_4$. Отсюда ${a_1a_2=a_3a_4}$.

Рассмотрим композицию инверсии и центральной симметрии с центром в $O$, оставляющую $\alpha$ на месте (такое преобразование существует т.~к. $O$ лежит внутри любой окружности I типа). В~этом случае любое решение переходит в противоположный сектор. Тогда $\omega_1$ меняется с $\omega_3$, а $\omega_2$ меняется с $\omega_4$. Отсюда $a_1a_3=a_2a_4$.

Из полученных равенств следует, что $a_1=a_4$ и $a_2=a_3$, что и требовалось.\qed

\begin{figure}[h!]
\begin{minipage}{0.4\linewidth}
\center{\includegraphics[width=1\linewidth]{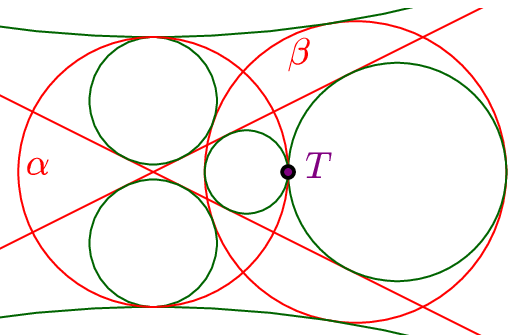}} \\ До инверсии в $T$
\end{minipage}
\hfill
\begin{minipage}{0.5\linewidth}
\center{\includegraphics[width=1\linewidth]{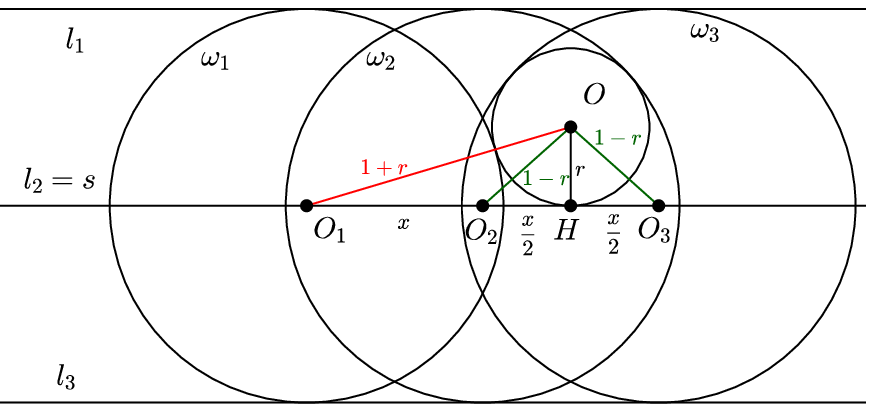}} \\ После инверсии в $T$
\end{minipage}
\caption{Одна окружность I типа и одна окружность II типа}
\label{ris:inv}
\end{figure}

Из этой леммы (после перенумерации секторов) получаем, что $\alpha$ и $\beta$ переходят в себя при симметрии относительно биссектрисы секторов I и III. Поэтому оба решения в I секторе касаются $\alpha$ в одной и той же точке $T$ (см. Рис.~\ref{ris:inv} слева). Заметим, что по лемме~\ref{lem:inter-4} $\alpha$ и $\beta$ имеют общие точки в секторах II и IV. Поэтому в секторе I они общих точек не имеют. Следовательно, $T\notin\beta$, и после инверсии в $T$ решения в I секторе и $\alpha$ переходят в 3 параллельные прямые $l_1$, $l_2$, $l_3$ (пусть $l_2$ --- образ $\alpha$, положим расстояние между $l_1$ и $l_3$ равным 2), а $\beta$ переходит в допустимую окружность $\omega_2$. Пересекающиеся прямые также переходят в допустимые окружности $\omega_1$ и $\omega_3$, причем в силу симметрии их центры равноудалены от центра $\omega_2$ (см. Рис.~\ref{ris:inv} справа). Обозначим это расстояние через $x$, а центры $\omega_i$ через $O_i$ $(i=1,2,3)$.

Обозначим через $s$ прямую, равноудаленную от $l_1$ и $l_3$. Четыре оставшихся решения в секторах II и IV переходят в окружности, касающиеся $\omega_1$, $\omega_2$ и $\omega_3$ в комбинациях <<$+--$>> и <<$--+$>>. По следствию из леммы~\ref{lem:algebraic} решения, подходящие для из этих наборов знаков, симметричны относительно $s$ --- линии центров $\omega_i$. Поэтому радиусы решений в каждой паре равны, а их линии центров вертикальны. Поэтому они должны касаться $l_2$ в одной и той же точке, и $l_2$ совпадает с $s$.

Обозначим радиус любого из решений через $r$, а его центр через $O$. Имеем $OO_1=1+r$, $OO_2=OO_3=1-r$. Пусть $OH$ --- высота равнобедренного $\bigtriangleup O_1OO_2$. Имеем $O_2H=HO_3=x/2$, $OH=r$. По теореме Пифагора в $\bigtriangleup O_1HO$ и $\bigtriangleup O_2HO$:
$$
\begin{cases}
OH^2+HO_1^2=OO_1^2\\
OH^2+HO_2^2=OO_2^2
\end{cases}
\begin{cases}
\left(\frac{3x}{2}\right)^2=(1+r)^2-r^2\\
\left(\frac{x}{2}\right)^2=(1-r)^2-r^2
\end{cases}
\begin{cases}
\frac{9x^2}{4}=1+2r\\
\frac{x^2}{4}=1-2r
\end{cases}
$$
Отсюда $x=2/\sqrt{5}$, $r=2/5$, что соответствует случаю №3 нашей классификации.

\boldmath\textbf{3. Обе окружности $\alpha$ и $\beta$ имеют I тип.}\unboldmath\ Тогда в каждом секторе не более двух решений и возможны распределения 1-2-1-2 и 1-1-2-2 (распределение 2-2-2-0 невозможно, т.~к. иначе по лемме~\ref{lem:inter-1} $\alpha$ и $\beta$ имеют 3 общие точки). Разберем их по отдельности.

\textit{3.1. Распределение 1-2-1-2.} Тогда $\alpha$ и $\beta$ касаются решений в секторах II и IV в комбинации <<$+--+$>> (в том порядке, в котором они касаются любой из прямых). По следствию из леммы~\ref{lem:algebraic} окружности $\alpha$ и $\beta$ симметричны относительно биссектрисы секторов II и IV. Получаем случай №2 в нашей классификации.

\textit{3.2. Распределение 1-1-2-2.} Тогда по лемме~\ref{lem:inter-1} $\alpha$ и $\beta$ имеют общую точку в секторах I и II. Сделаем инверсию в одной из них. Заметим, что одна из прямых разделяет точки пересечения $\alpha$ и $\beta$, а другая --- нет (см. Рис.~\ref{ris:lastcase}). Поэтому после инверсии одна из прямых перейдет в окружность I типа, а другая перейдет в окружность II типа. Этот случай уже рассмотрен.
\paragraph{Второй случай: никакие две окружности или прямые из исходных четырех не пересекаются.}
\begin{lemma}\label{lem:noninter}
Даны 3 окружности и 6 их решений. Тогда можно выбрать 2 исходные окружности и такую комбинацию знаков для них, что из 6 решений найдутся хотя бы 4 подходящих для выбранной пары и комбинации знаков.
\end{lemma}
\textbf{Доказательство леммы~\ref{lem:noninter}.} Рассмотрим таблицу из 3-х строк и 6-и столбцов, где строки соответствуют окружностям, а столбцы --- решениям. На пересечении i-ой строки и j-го столбца стоит знак, соответствующий касанию i-ой окружности с j-м решением. Достаточно доказать, что найдутся 2 строки и 4 столбца такие, что 4 пары знаков в выбранных столбцах совпадают с точностью до замены обоих знаков в паре на противоположные.

Инвертируем столбцы так, чтобы в 1-ой строке все знаки были <<$+$>>. Если тогда во 2-й или 3-й строке есть 4 одинаковых знака, то эта строка вместе с 1-ой является искомой. В~противном случае в каждой из них 3 знака <<$+$>> и 3 знака <<$-$>>. Пусть во 2-й строке в первых 3-х столбцах стоят <<$+$>>, а в остальных <<$-$>>. Если тогда в 3-й строке в первых 3-х столбцах одинаковые знаки, то 2-я и 3-я строки искомые. Пусть в 3-й строке в первых 3-х столбцах два одинаковых знака, скажем, <<$+$>>. Тогда в последних 3-х столбцах есть два знака <<$-$>> и эти 4 столбца и последние две строки являются искомыми --- лемма доказана.\qed
\begin{corollary}
Пусть в условиях леммы исходные окружности попарно не пересекаются. Тогда можно перевести некоторые две окружности в концентрические так, чтобы нашлись либо 4 решения типа A, либо 4 решения типа B.
\end{corollary}
\textbf{Доказательство следствия.} Всего возможно 2 комбинации из двух знаков с точностью до замены обоих на противоположные. Нетрудно видеть, что если исходные окружности не пересекаются, то одна из этих комбинаций соответствует тому, что решение разделяет эти окружности, а вторая --- что не разделяет. Поэтому переведя окружности из леммы в концентрические, мы получим 4 решения, которые либо все одновременно разделяют концентрические, либо нет, т.~е. принадлежат к типу A или B.\qed

Продолжим разбор второго случая теоремы~3. Применим следствие из леммы~\ref{lem:noninter} и переведем соответствующие окружности в концентрические. Остаются две исходные окружности $\alpha$ и $\beta$, порождающие 6 общих решений, из которых 4 имеют один и тот же тип. По лемме~\ref{lem:alternative} окружности $\alpha$ и $\beta$ либо лежат на одном диаметре, либо они лежат на перпендикулярных диаметрах. Разберем эти случаи по отдельности.

\boldmath\textbf{1. Окружности $\alpha$ и $\beta$ лежат на одном диаметре.}\unboldmath\ Как и в теореме~1, делаем инверсию в точке пересечения симметричных решений, переводя их в прямые $a$ и $b$. Поскольку имеется 6 решений, то по лемме~\ref{lem:twosectors} 4 исходные окружности должны иметь распределение 2-0-2-0 относительно $a$ и $b$, причем есть как решение I типа, так и решение II типа. После применения леммы~\ref{lem:beaty} имеем случай №1 нашей классификации.

\boldmath\textbf{2. Окружности $\alpha$ и $\beta$ лежат на перпендикулярных диаметрах.}\unboldmath\ Тогда 4 общих решения для $\alpha$ и $\beta$ имеют тип B, а их центры образуют прямоугольник. Также $\alpha$ и $\beta$ порождают еще 2 общих решения типа A.

\begin{figure}
\center{\includegraphics{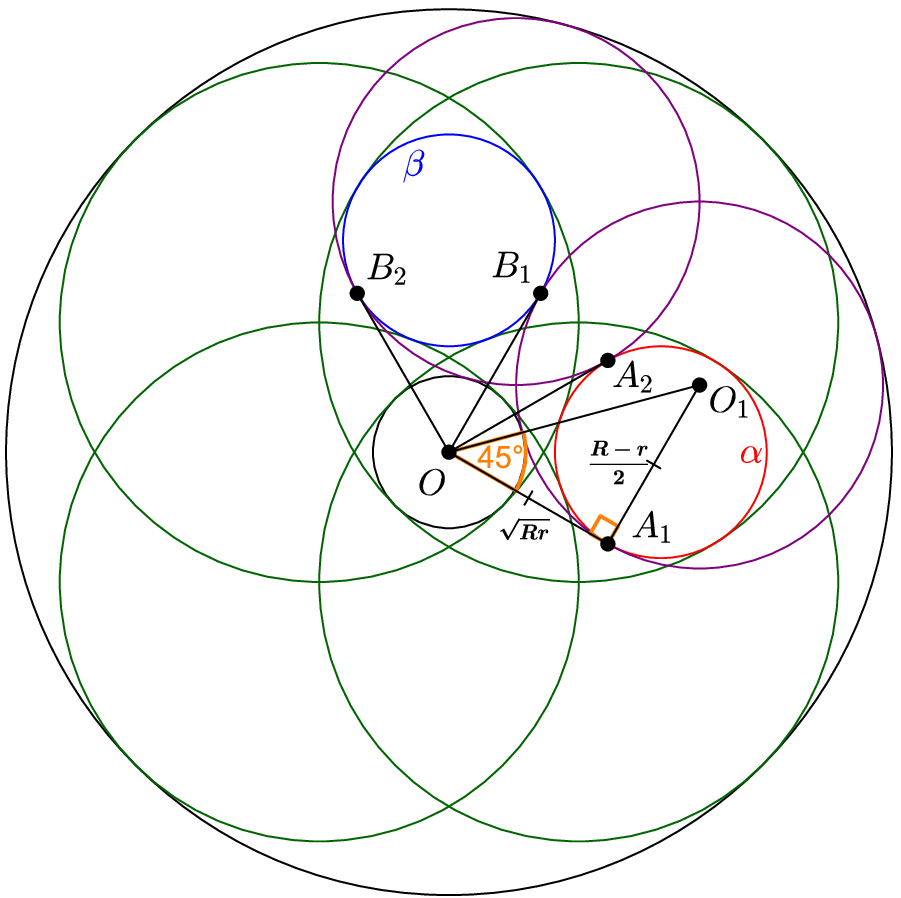}}
\caption{Случай непересекающихся исходных окружностей}
\label{ris:ConcExProof}
\end{figure}

Пусть $O$ --- центр концентрических окружностей, $R$ и $r$ --- их радиусы, $R>r$. Далее, пусть $A_1$, $B_1$ --- точки касания одного из решений типа A с $\alpha$ и $\beta$ соответственно, $A_2$, $B_2$ --- точки касания другого решения с $\alpha$ и $\beta$ соответственно.

Сделаем инверсию с центром в $O$ и радиусом $\sqrt{Rr}$. При такой инверсии все решения типа B перейдут в противоположные себе (т.~к. степень точки $O$ относительно всех решений типа B равна $-Rr$). Образы окружностей $\alpha$ и $\beta$ должны их касаться, поэтому $\alpha$ и $\beta$ перейдут в себя. Это значит, что степень точки $O$ относительно $\alpha$ и $\beta$ равна $Rr$. Однако степень точки $O$ относительно всех окружностей типа A также равна $Rr$, поэтому $OA_1$, $OA_2$ --- общие касательные к $\alpha$ и решениям типа A, а $OB_1$, $OB_2$ --- общие касательные к $\beta$ и решениям типа A, причем $OA_1=OA_2=OB_1=OB_2=\sqrt{Rr}$. Более того, $\angle A_1OB_1=\angle A_2OB_2$, т.~к. это угол, под которым решения типа A видны из $O$. Поэтому и $\angle A_1OA_2=\angle B_1OB_2$, а так как $OA_1=OB_1$, то радиусы $\alpha$ и $\beta$ равны. Тогда при повороте на $90^\circ$ относительно $O$ $\alpha$ и $\beta$ совпадут, а тогда $A_1$ и $B_1$ совпадут, т.~е. $\angle A_1OB_1=90^\circ$ (см. Рис.~\ref{ris:ConcExProof}).

Далее, пусть $O_1$ --- центр решения типа A. Он лежит на биссектрисе прямого угла $\angle A_1OB_1$. Поэтому $\bigtriangleup OA_1O_1$ --- прямоугольный равнобедренный, откуда
\begin{gather*}
OA_1=A_1O_1\\
\sqrt{Rr}=\frac{R-r}{2}\\
4Rr=(R-r)^2\\
R^2-6Rr+r^2=0
\end{gather*}
Следовательно, $R/r=3+2\sqrt{2}$ --- это случай №4 нашей классификации. Теорема~3 доказана.\qed

\section{Приложение: дальнейшие обобщения}
Что будет, если продолжить увеличение числа исходных окружностей? На самом деле ответить на этот вопрос можно исходя лишь из доказанных теорем. Действительно, пусть, например, дано 6 обобщенных окружностей. Если тогда существует 5 или более решений, то рассмотрим решения в качестве исходных окружностей, получая противоречие с теоремой~2. Ответ на вопрос о максимальном возможном числе решений для любого числа исходных окружностей содержится в Табл.~\ref{tab:bipart}.

Возможно ли тогда дальнейшее обобщение задачи Аполлония? Далее мы рассмотрим несколько путей.

\subsection{Обобщение задачи Штейнера}
Одним из обобщений задачи Аполлония является следующая
\begin{problem}[Штейнер]
На плоскости дано 5 коник, причем количество коник, касающихся каждой из данных, конечно. Какое максимально возможное значение может принимать это количество?
\end{problem}
Прямая аналогия заключается в том, что окружность задается тремя параметрами, а коника --- пятью. Поэтому получающееся здесь количество (3264, см.~\cite{bib:kiritch}) соответствует количеству решений классической задачи Аполлония (8). В связи с этим, возникает естественное аналогичное обобщение.
\begin{problem}[обобщенная задача Штейнера]
На плоскости дано $k\ge6$ коник, причем количество коник, касающихся каждой из данных, конечно. Какое максимально возможное значение может принимать это количество?
\end{problem}

\subsection{Графы касающихся окружностей}
Пусть $G$ --- граф без петель и кратных ребер. Будем говорить, что \emph{$G$ реализуем в касающихся окружностях}, если существует такой невырожденный набор окружностей на плоскости, что каждой вершине $G$ соответствует окружность, и окружности, соответствующие смежным вершинам, касаются.

Реализуемые графы --- очень интересные и недостаточно исследованные объекты. Вот почти все, что о них известно.
\begin{enumerate}
\item{Любой дистанционный граф реализуем.}
\item{Любой планарный граф реализуем (\emph{теорема Кобы-Андреева-Терстона об упаковке кругов}).}
\item{Граф $K_5$ нереализуем.}
\end{enumerate}
Напомним, что дистанционным называется граф, вершинами которого являются точки плоскости, а ребра проводятся между точками, находящимися на расстоянии 1.

Все утверждения, кроме второго, просты. Первое утверждение очевидно --- достаточно расположить окружности радиуса $1/2$ в вершинах графа. Третье утверждение доказывается от противного путем инверсии в точке касания любых двух окружностей. Эти окружности перейдут в параллельные прямые, а оставшиеся перейдут в три попарно касающиеся допустимые окружности, что очевидно невозможно.

\begin{table}[h!]
\begin{center}
\begin{tabular}{|c|c|c|c|c|c|c|c|c|}
\hline
$m, n$ & $\le2$ & $3$ & $4$ & $5$ & $6$ & $7$ & $8$ & $\ge9$ \\
\hline
$\le2$ & $+$ & $+$ & $+$ & $+$ & $+$ & $+$ & $+$ & $+$ \\
\hline
$3$ & $+$ & $+$ & $+$ & $+$ & $+$ & $+$ & $+$ & $-$ \\
\hline
$4$ & $+$ & $+$ & $+$ & $+$ & $+$ & $-$ & $-$ & $-$ \\
\hline
$5$ & $+$ & $+$ & $+$ & $-$ & $-$ & $-$ & $-$ & $-$ \\
\hline
$6$ & $+$ & $+$ & $+$ & $-$ & $-$ & $-$ & $-$ & $-$ \\
\hline
$7$ & $+$ & $+$ & $-$ & $-$ & $-$ & $-$ & $-$ & $-$ \\
\hline
$8$ & $+$ & $+$ & $-$ & $-$ & $-$ & $-$ & $-$ & $-$ \\
\hline
$\ge9$ & $+$ & $-$ & $-$ & $-$ & $-$ & $-$ & $-$ & $-$ \\
\hline
\end{tabular}
\end{center}
\caption{Реализуемость двудольных графов $K_{m, n}$}
\label{tab:bipart}
\end{table}

В~данной работе фактически были описаны все реализуемые двудольные графы. Результаты собраны в Табл.~\ref{tab:bipart}.

Дальнейшее исследование графов из окружностей также может быть рассмотрено как обобщение задачи Аполлония. Например, неизвестно, ограничено ли хроматическое число графов касающихся окружностей (это т. н. \emph{проблема Рингела}, см.~\cite{bib:ringel}).

\subsection{Высшие размерности}
Рассмотрение высших размерностей всегда является одним из естественных обобщений. Однако вопрос о реализуемости двудольных графов в касающихся сферах не представляет большого интереса ввиду следующего факта.
\begin{ntheorem}
Граф $K_{m, n}$ реализуем в касающихся сферах при любых $m, n$.
\end{ntheorem}
\textbf{Доказательство.} Достаточно рассмотреть $n$ сфер, вписанных в произвольный тор, и $m$ сфер, касающихся этого тора, с центрами на его оси.\qed

Данная теорема была предложена в качестве задачи на региональном этапе Всероссийской олимпиады школьников по математике 2016 (задача 11.6, см.~\cite{bib:vseros}).

Однако графы из касающихся сфер также могут быть исследованы. Например, для них актуален аналог проблемы Рингела.

Заслуживают внимания также обобщения задачи Аполлония в неевклидовых геометриях. 
\section*{Благодарности}
Я~признателен Ф.\,К.\,Нилову за постановку задачи и постоянное внимание к работе. Неоценимую помощь в подготовке данного текста оказал М.\,Б.\,Скопенков. Также хотелось бы поблагодарить М.\,А.\,Волчкевича за формулировку задачи, из развития которой впоследствии появилась данная работа.

\end{document}